\documentclass[11pt]{book}

\usepackage{amsfonts, amsmath, amssymb}  
\usepackage{times}

\usepackage{graphicx}
\DeclareGraphicsExtensions{.pdf,.png,.jpg}

\usepackage{amsfonts}
\usepackage{amssymb}
\usepackage{xfrac}
\usepackage{algorithm}
\usepackage[noend]{algpseudocode}
\usepackage{amsmath}
\usepackage{amsthm}
\usepackage{changepage}
\usepackage{lipsum}
\usepackage[export]{adjustbox}

\algrenewcommand\algorithmicrequire{\textbf{Precondition:}}
\algrenewcommand\algorithmicensure{\textbf{Postcondition:}}

\usepackage[T1]{fontenc}
\usepackage[utf8]{inputenc}
\usepackage{lmodern}
\usepackage[miktex]{gnuplottex}
\usepackage{graphicx}

\newtheorem*{theorem1}{Theorem 3.2}
\newtheorem*{lemma321}{Lemma 3.2.1}
\newtheorem*{lemma322}{Lemma 3.2.2}
\newtheorem*{lemma323}{Lemma 3.2.3}
\newtheorem*{lemma324}{Lemma 3.2.4}
\newtheorem*{lemma325}{Lemma 3.2.5}

\newtheorem*{theorem331}{Theorem 3.3.1}
\newtheorem*{propA}{(A)}
\newtheorem*{propB}{(B)}
\newtheorem*{propC}{(C)}
\newtheorem*{propD}{(D)}
\newtheorem*{propE}{(E)}

\newtheorem*{lemma332}{Lemma 3.3.1}
\newtheorem*{theorem332}{Theorem 3.3.1.1}
\newtheorem*{lemma3321}{Lemma 3.3.2}
\newtheorem*{lemma3322}{Lemma 3.3.3}
\newtheorem*{def1}{Step}
\newtheorem*{def2}{Position (of a kangaroo)}
\newtheorem*{def3}{Distance (between kangaroos)}

\newtheorem*{lemma3231}{Lemma 3.2.3.1}
\newtheorem*{lemma3232}{Lemma 3.2.3.2}

\begin{document}
\pagestyle{empty}

%: ----------------------------------------------------------------------
%:                  TITLE PAGE: name, degree,..
% ----------------------------------------------------------------------

\begin{center}

\vspace{1cm}

%%% Type the thesis title below%%%%%%%%%%%%%%%%
{\Huge         Kangaroo Methods for Solving the Interval Discrete Logarithm Problem}

\vspace{35mm} 

\includegraphics[width=2cm]{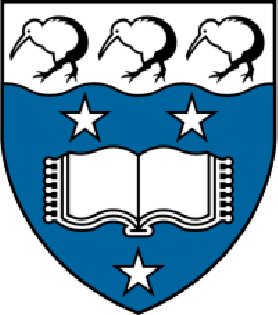}

 \vspace{45mm}

%%%%%Type Your Name Below%%%%%%%%%%%%
{\Large       Alex Fowler}

	\vspace{1ex}

Department of Mathematics

The University of Auckland

	\vspace{3ex}

 %%%%%Typing Your Supervisors Name Below%%%%%%%%%%%%
Supervisor:             Steven Galbraith

	\vspace{5mm}

A dissertation  submitted in partial fulfillment of the requirements for the degree of BSc(Hons)  in Applied Mathematics, The University of Auckland, 2014.

\end{center}

  \newpage

%: --------------------------------------------------------------
%:                  FRONT MATTER:  abstract,..
% --------------------------------------------------------------
\chapter*{Abstract}       
\setcounter{page}{1}
\pagestyle{headings}
% \pagenumbering{roman}

\addcontentsline{toc}{chapter}{Abstract}

The interval discrete logarithm problem is defined as follows: Given some $g,h$ in a group $G$, and some $N \in \mathbb{N}$ such that $g^z=h$ for some $z$ where $0 \leq z < N$, find $z$.\\ At the moment, kangaroo methods are the best low memory algorithm to solve the interval discrete logarithm problem. The fastest non parallelised kangaroo methods to solve this problem are the three kangaroo method, and the four kangaroo method. These respectively have expected average running times of $\big(1.818+o(1)\big)\sqrt{N}$, and $\big(1.714 + o(1)\big)\sqrt{N}$ group operations.\\It is currently an open question as to whether it is possible to improve kangaroo methods by using more than four kangaroos. Before this dissertation, the fastest kangaroo method that used more than four kangaroos required at least $2\sqrt{N}$ group operations to solve the interval discrete logarithm problem. In this thesis, I improve the running time of methods that use more than four kangaroos significantly, and almost beat the fastest kangaroo algorithm, by presenting a seven kangaroo method with an expected average running time of $\big(1.7195 + o(1)\big)\sqrt{N} \pm O(1)$ group operations. The question, 'Are five kangaroos worse than three?' is also answered in this thesis, as I propose a five kangaroo algorithm that requires on average $\big(1.737+o(1)\big)\sqrt{N}$ group operations to solve the interval discrete logarithm problem. 

%: --------------------------------------------------------------
%:                  END:  abstract
% --------------------------------------------------------------

%: ----------------------- contents ------------------------
\setcounter{secnumdepth}{3} % organisational level that receives a numbers
\setcounter{tocdepth}{3}    % print table of contents for level 3
\tableofcontents            % print the table of contents
% levels are: 0 - chapter, 1 - section, 2 - subsection, 3 - subsection

%: --------------------------------------------------------------
%:                  MAIN DOCUMENT SECTION
% --------------------------------------------------------------
	
\chapter{Introduction}%    \chapter{}  = level 1, top level
\section{Introduction} 
The interval discrete logarithm problem (IDLP) is defined in the following manner: Given a group $G$, some $g,h \in G$, and some $N \in \mathbb{N}$ such that $g^z = h$ for some $0 \leq z < N$, find $z$. In practice, $N$ will normally be much smaller than $\big|G\big|$, and $g$ will be a generator of $G$. The probability that $z$ takes any integer value between $0$ and $N-1$ is deemed to be uniform.\\
To our current knowledge, the IDLP is hard over some groups. Examples of such groups include Elliptic curves of large prime order, and $\mathbb{Z}_{p}^{*}$, where $p$ is a large prime. Hence, cryptosystems such as the Boneh-Goh-Nissim homomorphic encryption scheme [1] derive their security from the hardness of the IDLP. The IDLP also arises in a wide range of other contexts. Examples include, counting points on elliptic curves \cite{4}, small subgroup and side channel attacks [5,6], and the discrete logarithm problem with c-bit exponents [4]. The IDLP is therefore regarded as a very important problem in contemporary cryptography.\\
Kangaroo methods are the best generic low storage algorithm to solve the IDLP. In this thesis, I examine serial kangaroo algorithms, although all serial kangaroo methods can be parallelised in a standard way, giving a speed up in running time in the process. Currently, the fastest kangaroo algorithm is the 4 kangaroo method of Galbraith, Pollard and Ruprai [3]. On an interval of size $N$, this algorithm has an estimated average running time of $(1.714 + o(1))\sqrt{N}$ group operations and requires $O(\log(N))$ memory. It should be noted that over some specific groups, there are better algorithms to solve the IDLP. For example, over groups where inversion is fast, such as elliptic curves, an algorithm proposed by Galbraith and Ruprai in [2] has an expected average case running time of $(1.36 + o(1))\sqrt{N}$ group operations, while requiring a constant amount of memory. This method is much slower than the four kangaroo method over groups where inversion is slow however. Baby-step giant-step algorithms, which were first proposed by Shanks in [10], are also faster than kangaroo methods. The fastest Baby-step Giant-step algorithm was illustrated by Pollard in [8], and requires on average $4/3 \sqrt{N}$ group operations to solve the IDLP. However, these algorithms are unusable over large intervals, since they have $O(\sqrt{N})$ memory requirements. If one is solving the IDLP in an arbitrary group, on an interval of size $N$, one would typically use baby-step giant-step algorithms if $N < 2^{30}$, and would use the four kangaroo method if $N>2^{30}$.\\ 
The first kangaroo method was proposed by Pollard in 1978 in [9]. This had an estimated average running time of $3.3\sqrt{N}$ group operations [11]. The next improvement came from van Oorshot and Wiener in [12,13], with the introduction of an algorithm with an estimated average running time of $\big(2+o(1)\big)\sqrt{N}$ group operations. This was the fastest kangaroo method for over 15 years, until Galbraith, Pollard and Ruprai published their three and four kangaroo methods in [3]. These respectively require on average $\big(1.818+o(1)\big)\sqrt{N}$, and $\big(1.714+o(1)\big)\sqrt{N}$ group operations to solve the IDLP. The fastest kangaroo method that uses more than four kangaroos is a five kangaroo method, proposed in [3]. This requires at least $2\sqrt{N}$ group operations to solve the IDLP.\\It is currently believed, but not proven, that the four kangaroo method is the optimal kangaroo method. This is a major gap in our knowledge of kangaroo methods, and so the main question that this dissertation attempts to answer is the following. \begin{itemize}
\item \textbf{Question 1:} Can we improve kangaroo methods by using more than four kangaroos? \end{itemize}
This dissertation also attempts to answer the following lesser, but also interesting problem. \begin{itemize}
\item \textbf{Question 2:} Are five kangaroos worse than three? \end{itemize}
In this report, I attempt to answer these questions, by investigating five kangaroo methods in detail. I then state how a five kangaroo method can be adapted to give a seven kangaroo method, giving an improvement in running time in the process.
	
\chapter[Kangaroo Methods]{Current State of Knowledge of Kangaroo Methods}

In this section I will give a brief overview of how kangaroo methods work, and of our current state of knowledge of kangaroo algorithms. 

\section{General intuition behind kangaroo methods}
The key idea behind all kangaroo methods known today, is that if we can express any $x \in G$ in two of the forms out of $g^{p}h$,$g^{q}$ or $g^{r}h^{-1}$, where $p,q,r \in \mathbb{N}$, then one can find $z$, and hence solve the IDLP. In all known kangaroo methods, we have a herd of kangaroos who randomly 'hop' around various elements of $G$. Eventually, 2 different kangaroos can be expected to land on the same group element. If we define the kangaroo's walks such that all elements of a kangaroo's walk are in one of the forms out of $g^{p}h$,$g^{q}$ or $g^{r}h^{-1}$, then when 2 kangaroos land on the same group element, the IDLP may be able to be solved.  

\section{van Oorshot and Weiner Method}
The van Oorshot and Weiner method of [12,13] was the fastest kangaroo method for over 15 years. In this method, there is one 'tame' kangaroo (labelled $T$), and one 'wild' kangaroo (labelled $W$). Letting $t_{i}$ and $w_{i}$ respectively denote the group elements $T$ and $W$ are at after $i$ 'jumps' of their walk, $T$ and $W$'s walks are defined recursively in the following manner.
%\begin{algorithm}
%	\caption{}
%	\begin{algorithmic}
%	 \State $t_{0} = g^{\frac{N}{2}}$ %\Comment We start $T$'s walk at the group element $g^{N/2}$
%	\State $w_{0} = h$
%	\end{algorithmic}
%\end{algorithm}

%{
%\setlength{\interspacetitleruled}{-.4pt}%
%\begin{algorithm}
%  $t_{0} = g^{\frac{N}{2}}$ \Comment $T$ starts his walk at the group element $g^{\frac{N}{2}}$
%  bbb\;
%  ccc\;
%\end{algorithm}
%}%

\begin{itemize}
\item How $T$'s walk is defined
\begin{itemize}
\item $T$ starts his walk at $t_{0} = g^{\frac{N}{2}}$. \footnote{Note that the method assumes $N$ is even, so $N/2 \in \mathbb{N}$}
\item $t_{i+1} = t_{i}g^{n}$, for some $n \in \mathbb{N}$. 
\end{itemize}
\item How $W$'s walk is defined
\begin{itemize}
\item $W$ starts his walk at $w_{0} = h$.
\item $w_{i+1} = w_{i}g^{m}$, for some $m \in \mathbb{N}$.
\end{itemize}
\end{itemize}

One should note also that the algorithm is arranged so that $T$ and $W$ jump alternately. Clearly, all elements of $T$'s walk are expressed in the form $g^{q}$, while all elements of $W$'s walk are expressed in the form $g^{p}h$, where $p,q \in \mathbb{N}$. Hence when $T$ and $W$ both visit the same group element, we will have $w_{i} = t_{j}$, for some $i,j \in \mathbb{N}$, from which we can obtain an equation of the form $g^{q} = g^{p}h$, for some $p,q \in \mathbb{N}$, which implies $z = q-p$. \\
Now if we structure the algorithm so that the amount each kangaroo jumps by at each step is dependent only on its current group element, then from the point where both kangaroos have visited a common group element onwards, both kangaroos walks will follow the same path. As a consequence of this, we can detect when $T$ and $W$ have visited the same group element, while only storing a small number of the elements of each kangaroo's walk. All kangaroo methods employ this same idea, and this is why kangaroo methods only require $O(\log(N))$ memory. \\
To arrange the algorithm so that the amount each kangaroo jumps by at each step is dependent only on its current group element, we create a hash function $H$, which randomly assigns a 'step size' to each element of $G$. When a kangaroo lands on some $x \in G$, the amount it jumps forward by at that step is $H(x)$ (so its current group element is multiplied by the precomputed value of $g^{H(x)}$). %$g^{H(x)}$ will be precomputed, so each kangaroo's jump only requires a single group operation. If $g^{H(x)}$ wasn't precomputed, we would have to compute $g^{H(x)}$ at each step, which requires several multiplication operations.\\ 
To use this property so that the algorithm requires only $O(\log(N))$ memory, we create a set of 'distinguished points', $D$, where $D \subset G$. $D$ is defined such that $\forall x \in G$, the probability that $x \in D$ is $c \log(N)/\sqrt{N}$, for some $c > 0$. If a kangaroo lands on some $x \in D$, we first check to see if the other kangaroo has landed on $x$ also. If it has, we can solve the IDLP. If the other kangaroo hasn't landed on $x$, then if $T$ is the kangaroo that landed on $x$, we store $x$, the $q$ such that $x = g^{q}$, and a flag indicating that $T$ landed on $x$. On the other hand, if $W$ landed on $x$, we store $x$, the $p$ such that $g^{p}h = x$, and a flag indicating that $W$ landed on $x$. Hence we can only detect a collision after the kangaroos have visited the same distinguished point. At this stage, the IDLP can be solved. A diagram of the process the algorithm undertakes in solving the IDLP is shown below.\\

\includegraphics[width=140mm,height = 50mm]{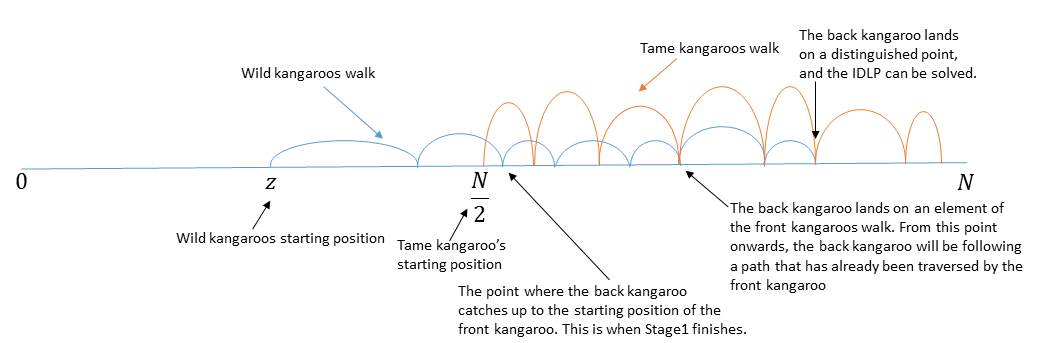}
%\begin{center} \texttt{van Oorshot and Weiner Method Diagram} \end{center}

To analyse the running time of the van Oorshot and Wiener method, we break the algorithm into the three disjoint stages shown in Stage 1, Stage 2, and Stage 3 below. In this analysis (and for the remainder of this thesis), I will use the following definitions.
\begin{def1} A period where each kangaroo makes exactly one jump.
\end{def1}
\begin{def2} A kangaroo is at position $p$ if and only if it's current group element is $g^p$. For instance, $T$ starts his walk at position $N/2$.
\end{def2}
\begin{def3} The difference between two kangaroos positions.
\end{def3}
We analyse the running time of the algorithm (and of all kangaroo methods) by considering the expected average number of group operations it requires to solve the IDLP. The reason for analysing the running time using this metric is explained in section 2.3.

\begin{itemize}
\item \textbf{Stage 1.} The period between when the kangaroos start their walks, and when the back kangaroo $B$ catches up to the front kangaroo $F$'s starting position. Now since the probability that $z$ takes any integer value $\big[0,N\big)$ is uniform, the average distance between $B$ and $F$ before they start their walks is $N/4$. Therefore, the expected number of steps required in stage 1 is $N/4m$, where $m$ is the average step size of the kangaroos walks
\item \textbf{Stage 2.} This is the period between when stage 1 finishes, and $B$ lands on a group element that has been visited by $F$. John Pollard showed experimentally in [8] that the number of steps required in this stage is $m$. One can see this intuitively in the following way. Once $B$ has caught up to $F$'s starting position, he will be jumping over a region in which $F$ has visited on average $1/m$ group elements. Hence the probability that $F$ lands on a element of $F$'s path at each step in Stage 2 is $1/m$. Therefore, we can expect $B$'s walk to join with $F$'s after $m$ steps.
\item \textbf{Stage 3.} This is the period between when stage 2 finishes, and $B$ lands on a distinguished point. Since the probability that any $x \in G$ is distinguished is $c\log(N)/\sqrt{N}$, for some $c > 0$, we can expect $B$ to make $\frac{1}{c\log(N)/\sqrt{N}}$ $= \sqrt{N}/c\log(N)$ steps in this stage.
\end{itemize}
%Now since the probability that $z$ takes any value $\big[0,N\big]$ is uniform, the average distance between $B$ and $F$ before they start their walks is $N/4$. Therefore, the expected number of steps required in stage 1 is $N/4m$, where $m$ is the average step size of the kangaroos walks. In \textbf{[Pollard Monopolies paper]}, John Pollard showed experimentally that the expected number of steps required in Stage 2 is $m$. One can see this intuitively in the following way.  Once $B$ has caught up to $F$'s starting position, he will be jumping over a region in which $F$ has visited on average $1/m$ group elements. Hence a heuristic estimate for the probability that $F$ lands on a element of $F$'s path at each step in Stage 2 is $1/m$. Therefore, we can expect $B$'s walk to join with $F$'s after $m$ steps. Once stage 2 has finished, since the probability that any $x \in G$ is distinguished is $c\log(N)/\sqrt{N}$, for some $c > 0$, we can expect $\sqrt{N}/c\log(N)$ steps to be made before $B$ lands on a distinguished point.
Hence we can expect the algorithm to require $N/4m + m + \sqrt{N}/c\log(N)$ steps to solve the IDLP. Now since each of the two kangaroos make one jump at each step, and in each jump a kangaroo makes we multiply two already known group elements together, each step requires two group operations. Hence the algorithm requires $2\big(N/4m + m + \sqrt{N}/c\log(N)\big)$ group operations to solve the IDLP. The optimal choice of $m$ in this expression is $m = \sqrt{N}/2$, which gives an expected average running time of $\big(2 + 1/c\log(N)\big)\sqrt{N}$ $= \big(2+o(1)\big)\sqrt{N}$ group operations.\\ Note that this analysis (and the analysis of the three and four kangaroo methods) ignores the number of group operations required to initialise the algorithm (in the initialisation phase we find the starting positions of the kangaroos, assign a step size to each $x \in G$, and precompute the group elements $g^{H(x)}$ for each $x \in G$). In section 3.3, I show however that the number of group operations required in this stage is constant. 
 
\section{How we analyse the running time of Kangaroo Methods} I can now explain why we analyse running time of kangaroo methods in terms of the expected average number of group operations they require to solve the IDLP. \\
\textbf{Why only consider group operations?} Generally speaking, in kangaroo methods, the computational operations that need to be carried out are group operations (e.g. multiplying a kangaroo's current group element to move it around the group), hashing (e.g. computing the step size assigned to a group element), and memory access comparisons (e.g. when a kangaroo lands on a distinguished point, checking to see if that distinguished point has already been visited by another kangaroo). The groups which one is required to solve the IDLP over are normally elliptic curve groups of large prime order, or $\mathbb{Z}_{p}^{*}$, for a large prime $p$ [11]. Group operations over such groups require far more computational operations than hashing, or memory access comparisons do. Hence when analysing the running time, we get an accurate approximation to the number of computational operations required by only counting the number of group operations. Counting the number of hashing operations, and memory access comparisons increases the difficulty of analysing the running time substantially, so it is desirable to only count group operations. \\ \textbf{Why consider the \emph{expected average} number of group operations?}
In the van Oorshot and Weiner method (and in all other kangaroo methods), the hash function which assigns a step size to each element of $G$ is chosen randomly. Now for fixed $z$ and varied hash functions, there are many possibilities for how the walks of the kangaroos will pan out. Hence in practise, the number of group operations until the IDLP is solved can take many possible values for each $z$. Hence we consider the \emph{expected} running time for each $z$, as being the average running time across all possible walks the kangaroos can make for each $z$. \\ Now $z$ can take any integer value between $0$ and $N-1$ with equal probability. Hence we refer to the expected \emph{average} running time, as being the average of the expected running times across all $z \in \mathbb{N}$ in the interval $\big[0,N)$.

\section{Three Kangaroo Method}
The next major breakthrough in kangaroo methods came from Galbraith, Pollard and Ruprai in [3] with the introduction of their three kangaroo method. The algorithm assumes that $g$ has odd order, and that $10 | N$. The method uses three different types of kangaroos, labelled $W_{1},W_{2}$ and $T$. All of the elements of $W_{1}$,$W_{2}$ and $T$'s walks are respectively expressed in the forms $g^{p}h$, $g^{r}h^{-1}$, and $g^{q}$, where $p,q,r \in \mathbb{N}$. If any pair of kangaroos collides, then we can express some $x \in G$ in 2 of the forms out of $g^{p}h$,$g^{q}$ and $g^{r}h^{-1}$. If we have $x = g^{p}h = g^{q}$, then $z = p-q$, while if we have $x = g^{q} = g^{r}h^{-1}$, then $z = r-q$. On the other hand, if we have $x = g^{p}h = g^{q}h^{-1}$, then since $g$ has odd order, we can solve $z = 2^{-1}\big(q-p\big) \mod(\big|g\big|)$. Hence a collision between any pair of kangaroos can solve the IDLP. 
\\To enable us to express the elements of the kangaroos walks in these ways, as in the van Oorshot and Wiener method, we create a hash function $H$ which assigns a step size to each element of $G$. Defining $W_{1,i}$,$W_{2,i}$ and $T_{i}$ to respectively denote the group element $W_{1}$, $W_{2}$ and $T$ are at after $i$ jumps in their walk, we then define the walks of the kangaroos in the following way.\begin{itemize}
\item How $W_{1}$'s walk is defined
\begin{itemize}
\item $W_{1}$ starts at the group element $W_{1,0} = g^{-N/2}h$
\item To move $W_{1}$ to the next step, $W_{1,i+1} = W_{1,i}g^{H(W_{1,i+1})}$ 
\end{itemize}
\item How $W_{2}$'s walk is defined
\begin{itemize}
\item We start $W_{2}$ at the group element $W_{2,0}=g^{N/2}h^{-1}$
\item To move $W_{2}$ to the next step, $W_{2,i+1} = W_{2,i}g^{H(W_{2,i})}$
\end{itemize}
\item How $T$'s walk is defined
\begin{itemize}
\item We start $T$ at $T_{0} = g^{3N/10}$.
\item To move $T$ to the next step, $T_{i+1} = T_{i}g^{H(T_{1})}$
\end{itemize}
\end{itemize}
As in the van Oorshot and Wiener method, the algorithm is arranged so that the kangaroos jump alternately. One can see that $W_{1}$,$W_{2}$, and $T$ start their walks respectively at the positions $z-N/2$,$N/2-z$, and $3N/10$. %maybe have diagram of three roo method that I used in the powerpoint here,
Hence if we let $d_{T,W_{1}}$,$d_{T,W_{2}}$, and $d_{W_{1},W_{2}}$ be the functions for the initial distances between $T$ and $W_{1}$, $T$ and $W_{2}$, and $W_{1}$ and $W_{2}$ over all $0 \leq z < N$,
 then $d_{T,W_{1}}(z) = \big|(z-N/2)-(3N/10)\big| = \big|z - 4N/5\big|$,
 $d_{T,W_{2}}(z) = \big|N/5-z\big|$, and $d_{W_{1},W_{2}}(z) = \big|2z-N\big|$. We also define $d_{C}$ to be the function which denotes the initial distance between the closest pair of kangaroos over all $0 \leq z < N$. Hence 
$d_{C}(z) = \min\big\{d_{T,W_{1}}(z),d_{T,W_{2}}(z),d_{W_{1},W_{2}}(z)\big\}$. The starting positions of the kangaroos in this method are chosen so that the average distance between the closest pair of kangaroos (the average of $d_{C}(z)$ for $0 \leq z < N$) is minimised. A diagram of the distance between all pairs of kangaroos, and of $d_{C}$ is shown below.\\

%\begin{figure}
\includegraphics[width = 12cm, height = 6cm]{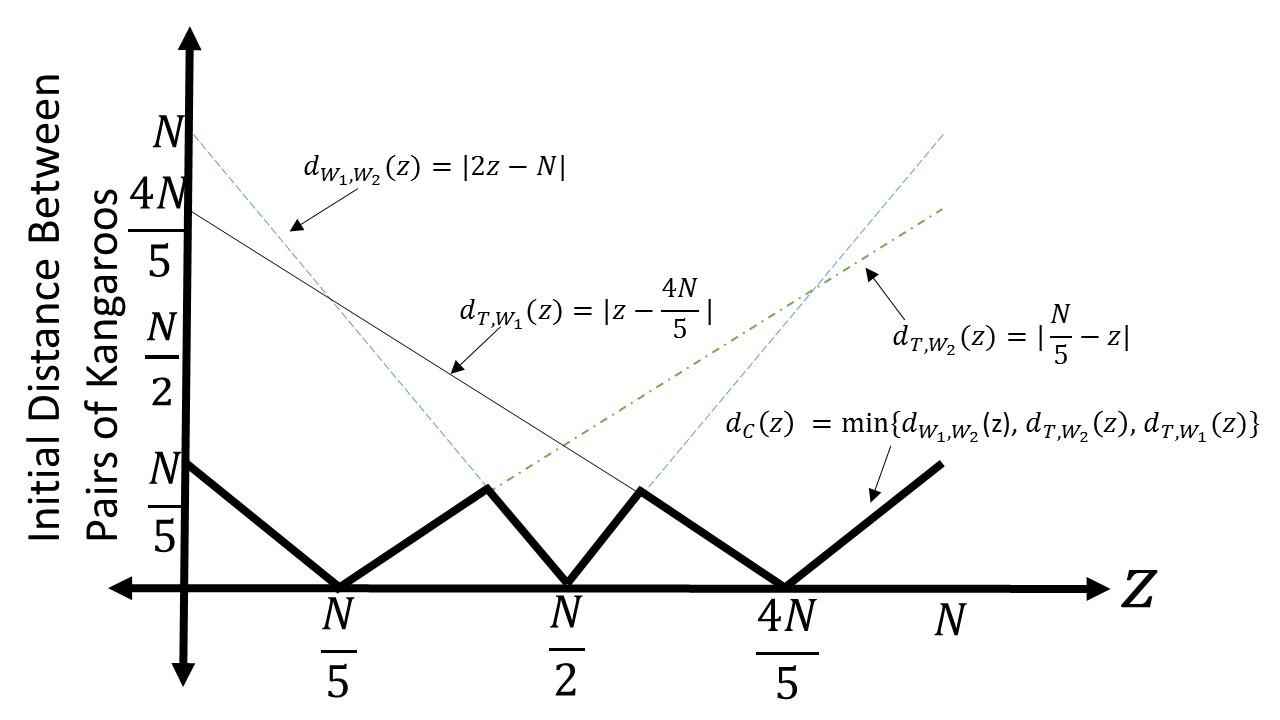}
%\caption{}\label{F1}
%\end{figure}

The following table shows the the formula for $d_{C}$, what $C$ (the closest pair of kangaroos) is, what $B$ (the back kangaroo in $C$) is, and what $F$ (the front kangaroo in $C$) is, across all $0 \leq z < N$.\\ 

\begin{tabular}[center]{|l|l|l|l|l|}
\hline
& $d_{C}(z)$ & $C$ & $B$ & F\\ \hline
$0 \leq z \leq N/5$ & $d_{T,W_{2}}(z) = N/5 - z$ & $T$ and $W_{2}$ & $T$ & $W_{2}$ \\ \hline
$N/5 \leq z \leq 2N/5$ & $d_{T,W_{2}}(z) = z - N/5$ & $T$ and $W_{2}$ & $W_{2}$ & $T$ \\ \hline
$2N/5 \leq z \leq N/2$ & $d_{W_{1},W_{2}} = N - 2z$ & $W_{1}$ and $W_{2}$ & $W_{1}$ & $W_{2}$ \\ \hline
$N/2 \leq z \leq 3N/5$ & $d_{W_{1},W_{2}} = 2z - N$ & $W_{1}$ and $W_{2}$ & $W_{2}$ & $W_{1}$ \\ \hline
$3N/5 \leq z \leq 4N/5$ & $d_{T,W_{1}}(z) = 4N/5 - z$ & $T$ and $W_{1}$ & $W_{1}$ & $T$ \\ \hline
$4N/5 \leq z \leq N/2$ & $d_{T,W_{1}}(z) = z - 4N/5$ & $T$ and $W_{1}$ & $T$ & $W_{1}$ \\ \hline
\end{tabular}

%One can see that for all $0 \leq z \leq 2N/5$, $W_{2}$ and $T$ is the closest pair (so $d(z)=d_{T,W_{2}}$), while for all $2N/5 \leq z \leq 3N/5$, $W_{1}$ and $W_{2}$ are the closest pair (so $d(z) = d_{W_{1},W_{2}}$), while $\forall \quad 3N/5 \leq z < N$, $T$ and $W_{1}$ are the closest pair (so $d(z) = d_{T,W_{1}}$). 
The expected number of steps until the IDLP is solved from a collision between the closest pair can be analysed in the same way as the running time was analysed in the van Oorshot and Weiner method. \begin{itemize}
\item \textbf{Stage 1.} The period between when the kangaroos start their walks, and when $B$ catches up with $F$'s starting position. The average of $d_{C}$ can easily be seen to be $N/10$. Hence if we let $m$ be the average step size used, the expected number of steps for $B$ to catch up to $F$'s starting position is $N/10m$.
\item \textbf{Stage 2.} The period between when stage 1 finishes, and when $B$ lands on an elements of $F$'s walk. The same analysis as was applied in the van Oorshot and Weiner method shows that the expected number of steps required in this stage is $m$.
\item \textbf{Stage 3.} The period between when $B$ lands on an element of $F$'s path, and $B$ lands on a distinguished point. In the three kangaroo method, the probability of a group element being distinguished is the same as it is in the van Oorshot and Wiener method, so the expected number of steps required in this stage is $\sqrt{N}/c\log(N)$.\end{itemize}

If we make the pessimistic assumption that the IDLP will always be solved from a collision between the closest pair of kangaroos, we can expect the algorithm to require $N/10m + m + \sqrt{N}/c\log(N)$ steps to solve the IDLP. This expression is minimised when $m$ is taken to be $\sqrt{N/10}$. In this case, the algorithm requires $\big(2\sqrt{1/10}+o(1))\big)\sqrt{N}$ steps to solve the IDLP. Since there are three kangaroos jumping at each step, the expected number of group operations until the IDLP is solved is $\big(1.897+o(1)\big)\sqrt{N}$ group operations.\\When Galbraith, Pollard and Ruprai considered the expected number of group operations until the IDLP was solved from a collision between any pair of kangaroos (so not just from a collision between the closest pair), they found through a complex analysis, that the three kangaroo method has an expected average case running time of $\big(1.818+o(1)\big)\sqrt{N}$ group operations. This is a huge improvement on the running time of the van Oorshot and Weiner method.  
%\includegraphics[width=125mm]{picture3.png}
%\begin{center}
%\texttt{Diagram of the Three Kangaroo Method}
%\end{center}

 %$W_{1}$, $W_{2}$ and $T$'s walks are respectively defined to start at the elements $g^{-N/2}h$, $g^{N/2}h^{-1}$, and $g^{3N/10}$. This corresponds to starting the kangaroos at the positions $z-N/2$,$N/2 - z$, and $3N/10$. The walks of the kangaroos are defined identically to how they were in the van Oorshot and Weiner method in every other way.

\section{Four Kangaroo Method}
The Four kangaroo method is a very simple, but clever extension of the 3 kangaroo method. As in the three kangaroo method, we start three kangaroos, $T_{1}$, $W_{1}$ and $W_{2}$ at the positions $\sfrac{3N}{10}$, $z-\sfrac{N}{2}$, and $\sfrac{N}{2}-z$ respectively. Here however, we add in one extra tame kangaroo ($T_{2}$), who starts his walk at $\sfrac{3N}{10} + 1$. One can see that the starting positions of $W_{1}$ and $W_{2}$ have the same parity, while exactly one of $T_{1}$ and $T_{2}$'s starting positions will have the same parity as $W_{1}$ and $W_{2}$'s starting positions. Therefore, if the step sizes are defined to be even, then in any walk, both of the wild, and one of the tame kangaroos will be able to collide, while one of the tame kangaroos will be unable to collide with any other kangaroo. Therefore, the three kangaroos that can collide are effectively simulating the three kangaroo method, except over an interval of half the size. Hence, from the analysis of the three kangaroo method, the three kangaroos that can collide in this method require $(1.818 + o(1))\sqrt{N/2}$ group operations to solve the IDLP. However, since there is one 'useless' kangaroo that requires just as many group operations as the three other useful kangaroos, the expected number of group operations required to solve the IDLP by the four kangaroo method is $\sfrac{(3+1)}{3}(1.818 + o(1))\sqrt{N/2} = (1.714 + o(1))\sqrt{N}$ group operations. 

\section{Can we do better by using more kangaroos?}

I now return to the main question this dissertation seeks to address. The answer to this question is not clear through intuition, since there are arguments both for and against using more kangaroos.\\On the one hand, if one uses more kangaroos, we both increase the number of pairs of kangaroos that can collide, and we can make the kangaroos closer together. Therefore, by using more kangaroos, the number of steps until the first collision occurs will decrease. However, by increasing the number of kangaroos, there will more kangaroos jumping at each step, so the number of group operations required at each step increases.

\chapter{Five Kangaroo Methods}

This section will attempt to answer the two main questions of this dissertation (see Question 1 and Question 2 in the introduction), by investigating kangaroo methods which use five kangaroos.

In this section, I will answer Question 2, and partially answer Question 1, by presenting a five kangaroo algorithm which requires on average $\big(1.737 + o(1)\big)\sqrt{N} \pm O(1)$ group operations to solve the IDLP. To find a five kangaroo algorithm with this running time, I answered the following questions, in the order stated below. \begin{itemize}

%The main questions to consider when creating a 5 kangaroo algorithm are the following.

\item How should the walks of the kangaroos be defined in a 5 kangaroo algorithm?
\item How many kangaroos of each type should be used?
\item Where abouts should the kangaroos start their walks?
\item What average step size should be used?
\end{itemize}

\section{How the Walks of the Kangaroos are Defined}

%In my 5 kangaroo algorithm, I will use the same three types of kangaroos as those that were used in the three and four kangaroo methods. There will be two different types of 'wild' kangaroo (called \textsc{wild1} and \textsc{wild2}), and one type of 'tame' kangaroo (labelled \textsc{tame}). Kangaroos of types \textsc{wild1} and \textsc{wild2} will start their walks at the group elements $g^{p}h$, and $g^{r}h^{-1}$ respectively, while \textsc{tame} kangaroos will start their walks at $g^{q}$, where $p,q,r \in \mathbb{N}$. This corresponds to staring the kangaroos at the positions $p+z$,$r-z$ and $q$. The kangaroos will then walk around the group, jumping one after the other, in some specified order. As before, a hash function $H$ will assign a step size to each group element in $G$, and the amount each kangaroo jumps by at any step will be $g$ to the power of the step size assigned to its current group element. Therefore, the elements of the walk of a kangaroo of types \textsc{tame}, \textsc{wild1} and \textsc{wild2} will be expressed in the forms $g^{n_{1}}$,$g^{n_{2}}h$, and $g^{n_{3}}h^{-1}$ respectively, where $n_{1}$,$n_{2}$,$n_{3} \in \mathbb{N}$. The mean step size will be defined to be $c_{m}\sqrt{N}$, for a constant $c_{m}$. 

I will first investigate 5 kangaroo methods where a kangaroo's walk can be defined in one of the same three ways as they were in the three and four kangaroo methods. This means, that a kangaroo can either be of type \textsc{Wild1},\textsc{Wild2}, or \textsc{Tame}, where the types of kangaroos are defined in the following way.\begin{itemize}
\item \textbf{Wild1 Kangaroo} - A kangaroo for which we express all elements of its walk in the form $g^{p}h$, where $p \in \mathbb{N}$.
\item \textbf{Wild2 Kangaroo} - A kangaroo for which we express all elements of its walk in the form $g^{r}h^{-1}$, where $r \in \mathbb{N}$.
\item \textbf{Tame Kangaroo} - A kangaroo for which we express all elements of its walk in the form $g^{q}$, where $q \in \mathbb{N}$.
\end{itemize}
As in the three and four kangaroo methods, there will be a hash function $H$ which assigns a step size to each $x \in G$. To walk the kangaroos around the group, if $x$ is a kangaroos current group element, the group element it will jump to next will be $xg^{H(x)}$. As in all previous kangaroo methods, the algorithm will be arranged so that the kangaroos each take one jump during each step of the algorithm.
 
\section{How many kangaroos of each type should be used?}

If we let $N_{W1}$, $N_{W2}$ and $N_{T}$ respectively denote the number of \textsc{wild1}, \textsc{wild2}, and \textsc{tame} kangaroos used in any 5 kangaroo method. Then $N_{W1} + N_{W2} + N_{T} = 5$. The following theorem will prove to be very useful in working out how many kangaroos of each type should be used, given this constraint.

\begin{theorem1} Let $A$ be any 5 kangaroo algorithm, where the kangaroos can be of type \textsc{Tame}, \textsc{Wild1}, or \textsc{Wild2}. Then the expected number of group operations until the closest 'useful' pair of kangaroos in $A$ collides is no less than $10\sqrt{\frac{N}{2N_{W1}N_{W2} + 4N_{T}(N_{W1}+N_{W2})}}$ (A pair is called useful if the IDLP can be solved when the pair collides). \end{theorem1}

My proof of this requires Lemma 3.2.1, Lemma 3.2.2, Lemma 3.2.3, Lemma 3.2.4, and Lemma 3.2.5. Before stating and proving these lemmas, I will make the following two important remarks. \begin{itemize}
\item \textbf{Remark 1.} I showed in my description of the three kangaroo method how a collision between kangaroos of different types could solve the IDLP. On the other hand, we can gain no information about $z$ from a collision between kangaroos of the same type. Hence a pair of kangaroos is 'useful' if and only if it features two kangaroos of different types.
\item \textbf{Remark 2.} In this section, for any choice of starting positions, I will let $d$ (where $d$ is a function of $z$) be the function that denotes the initial distance between the closest useful pair of kangaroos across all $0 \leq z < N$. $d$ is analogous to the function $d_{C}$ in the three kangaroo method, where for any $z$ with $0 \leq z < N$, $d(z)$ will be defined to be the smallest distance between pairs of kangaroos of different types, at that particular $z$. Note that $d(z)$ is completely determined by the starting positions of the kangaroos.
\end{itemize}

\begin{lemma321}For any choice of starting positions in a 5 kangaroo algorithm, the minimal expected number of group operations until the closest useful pair collides is $10 \sqrt{Ave(d(z))}$, where $Ave(d(z))$ is the average starting distance between the closest useful pair of kangaroos over all instances of the IDLP (i.e. over all $z \in \mathbb{N}$ where $z \in \big[0,N\big) = \big[0,N-1\big]$).\end{lemma321}

\begin{proof}[Proof of Lemma 3.2.1] Let $A$ be a 5 kangaroo algorithm that starts all kangaroos at some specified choice of starting positions, and let $m$ be the average step size used in $A$. Also let $d(z)$ be defined as in remark 2. For each $z$ with $0 \leq z < N$, using an argument very similar to that provided in section 2.2, the expected number of steps until the closest useful pair collides for this specific $z$ is $d(z)/m + m$. Since 5 kangaroos jump at each step, the expected number of group operations until the closest useful pair collides for this $z$ is $5(d(z)/m + m)$. Therefore, the expected average number of group operations until the closest useful pair collides across all instances of the IDLP (i.e. across all $z \in \mathbb{N}$ with $0 \leq z < N$) is $$\frac{1}{N} \sum\limits_{z=0}^{N-1} 5\left(\frac{d(z)}{m}+m\right)$$ $$\approx \frac{1}{N} \int_0^{N-1} 5\left(\frac{d(z)}{m} + m\right) dz$$ $$= 5\left(\frac{\frac{1}{N} \int_0^{N-1} d(z) dz}{m} + m\right)$$ $$\approx 5\left(\frac{\frac{1}{N}\sum\limits_{z=0}^{N-1} d(z)}{m} + m\right) = 5\left(\frac{Ave(d(z))}{m} + m\right).$$ Now simple differentiation shows that the $m$ that minimises this is $m = \sqrt{Ave(d(z))}$. Substituting in $m = \sqrt{Ave(d(z))}$ gives the required result. \end{proof}

\begin{lemma322} For each $z$, $d(z)$ is either of the form $|C_{1} \pm z|$ or $|C_{2} \pm 2z|$, where $C_{1}$ and $C_{2}$ are constants independent of $z$. Furthermore, letting $p_{g_{1}}$ and $p_{g_{2}}$ respectively denote the number of pairs with initial distance function of the form $|C \pm z|$, and $|C \pm 2z|$, $p_{g_{1}} = N_{T}(N_{W1} + N_{W2})$, and $p_{g_{2}} = N_{W1}N_{W2}$.\end{lemma322}

\begin{proof}[Proof of Lemma 3.2.2] Since a pair of kangaroos is useful if and only if it features two kangaroos of different types, a pair is useful if and only if it is a \textsc{tame}/\textsc{wild1}, a \textsc{tame}/\textsc{wild2}, or a \textsc{wild1}/\textsc{wild2} pair (here a \textsc{type1}/\textsc{type2} pair means a pair of kangaroos featuring one kangaroo of \textsc{type1}, and another kangaroo of \textsc{type2}). Since all \textsc{wild1},\textsc{wild2} and \textsc{tame} kangaroos start their walks respectively at group elements of the form $g^{p}h$,$g^{r}h^{-1}$ and $g^{q}$, for some $p,q,r \in \mathbb{N}$, the starting positions of all \textsc{Wild1},\textsc{Wild2} and \textsc{Tame} kangaroos will be of the forms $p + z$, $r - z$, and $q$ respectively. Hence the initial distance functions between \textsc{tame}/\textsc{wild1}, \textsc{tame}/\textsc{wild2}, and \textsc{wild1}/\textsc{wild2} pairs are respectively $|p + z - q|$, $|r - z + q|$, and $|p + z - (r - z)|$. Therefore, the distance function between all \textsc{tame}/\textsc{wild1} and \textsc{tame}/\textsc{wild2} pairs can be expressed in the form $|C_{1} \pm z|$, where $C_{1}$ is independent of $z$. The number of such pairs is $N_{T}N_{W1} + N_{T}N_{W2}$. On the other hand, the initial distance function between a \textsc{wild1}/\textsc{wild2} pair can be expressed in the form $|C_{2} \pm 2z|$, where $C_{2}$ is also independent of $z$. The number of such pairs is $N_{W1}N_{W2}$.
\end{proof}

%\begin{lemma324} The area under $d$ is minimised when each pair is the closest pair only over no more than one region. Furthermore, the area under $d$ is minimised when over each region for which a pair is closest over 
%\end{lemma324}

From the graph of $d_{C}$ in section 2.4, we can see that the function for the distance between the closest pair of kangaroos in the three kangaroo method is a sequence of triangles. The same is generally true in five kangaroo methods. In lemma 3.2.3, I will show that the area under the function $d$ is minimised in five kangaroo methods when $d$ is a sequence of triangles.

%THIS PREVIOUS PAPRAGRAPH REDONE:
%From the graph of $d_{C}$ in section 2.4, we can see that the function for the distance between the closest pair of kangaroos in the three kangaroo method is a sequence of triangles. The same is generally true in five kangaroo methods. In this section, I will make the assumption that the area under $d$ will be minimised when $d$ is a sequence of triangles. This is highly likely to be a valid assumption

\begin{lemma323} Suppose $d_{1,1},d_{1,2},d_{1,3},...,d_{p_{1}}$ and $d_{2,1},d_{2,2},...,d_{2,p_{2}}$ are functions of $z$, defined on the interval $\left[0,N\right)$, where for all $i$ and $j$, $d_{1,i}(z) = |C_{i} \pm z|$, and $d_{2,j}(z) = |C_{2} \pm 2z|$, where $C_{i}$ and $C_{j}$ are constants. Let $d(z)$ be defined such that $\forall \ 0 \leq z < N$, $d(z) = \min\left\{d_{1,1}(z),d_{1,2}(z),...,d_{1,p_{1}}(z),d_{2,1}(z),...,d_{2,p_{2}}(z)\right\}$. Then assuming that we have full control over the constants $C_{i}$ and $C_{j}$, in every case where $d$ is not a sequence of triangles, the area under $d$ can be decreased by making $d$ into a sequence of triangles. This can be done by changing some of the constants $C_{i}$ and $C_{j}$. \end{lemma323}

\begin{proof}
By $d$ being a sequence of triangles, I mean that if we shade in the region beneath $d$, then the shaded figure is a sequence of triangles (see figure 3.2(a)). Now $d$ is a sequence of triangles if and only if $S_{1}$,$S_{2}$ and $S_{3}$ hold, where $S_{1}$ is the statement '$d(0) = 0$, or $d(0) > 0$ and $d\ensuremath{'}(0) < 0$', $S_{2}$ is the statement '$d(N-1) = 0$, or $d(N-1)> 0$ and $d\ensuremath{'}(N-1) > 0$', and $S_{3}$ is the statement 'for every $z$ where the gradient of $d$ changes, the sign of the gradient of $d$ changes'. This is equivalent to the statement, 'for every $z$ where the function which is smallest changes (from say $d_{a_{1},b_{1}}$ to $d_{a_{2},b_{2}}$), the gradients of both $d_{a_{1},b_{1}}$ and $d_{a_{2},b_{2}}$ have opposite sign at $z$'.\\Hence if $d$ is not a sequence of triangles, $d$ must not satisfy at least one of the properties out of $S_{1}$ (see figure 3.2(b)),$S_{2}$ (see figure 3.2(c)), or $S_{3}$ (see figures 3.2 (d) and (e)). %'for every $z$ where the function which is smallest changes (from say $d_{a_{1},b_{1}}$ to $d_{a_{2},b_{2}}$), the gradients of both $d_{a_{1},b_{1}}$ and $d_{a_{2},b_{2}}$ have opposite sign at $z$'. Figure 3.2(b) is an example of when $S_{1}$,$S_{2}$ and $S_{3}$ all hold.\\ Now suppose $d$ is not a sequence of triangles. Then at least one of $S_{1}$,$S_{2}$ or $S_{3}$ doesn't hold.\\

%Suppose $d$ is not a sequence of triangles.  For $d$ to not be a sequence of triangles, then either $d(0) > 0$, and $d(0)\ensuremath{'} > 0$, (see figure 3.2(a)), $d(N) > 0$ and $d(N)\ensuremath{'} < 0$ (see figure 3.2(b)) or there are functions $d_{a_{1},b_{1}},d_{a_{2},b_{2}}$ out of  $d_{1,1},...,d_{1,p_{g_{1}}},d_{2,1},...,d_{2,p_{g_{2}}}$ for which at some $z$, $d$ changes from being equal to $d_{a_{1},b_{1}}$ to being equal to $d_{a_{2},b_{2}}$, and the gradients of $d_{a_{1},b_{1}}$ and $d_{a_{2},b_{2}}$ have the same sign (see figure 3.2(c) and figure 3.2(d)). If the converse of this holds (so $S_{1}$, $S_{2}$, and $S_{3}$ all hold, where $S_{1}$ is the statement '$d(0) = 0$ or ($d(0) > 0$ and $d(0)\ensuremath{'} < 0$)', $S_{2}$ is the statement '$d(N) = 0$ or ($d(N)> 0$ and $d(0)\ensuremath{'} < 0$)', and $S_{3}$ is the statement 'for every $z$ where the function which is smallest changes (from say $d_{a_{1},b_{1}}$ to $d_{a_{2},b_{2}}$), the gradients of both $d_{a_{1},b_{1}}$ and $d_{a_{2},b_{2}}$ have opposite sign at $z$'), then $d$ is a sequence of triangles (see figure 3.2(e)).\\ 
\includegraphics[height = 4.5cm, width = 12cm]{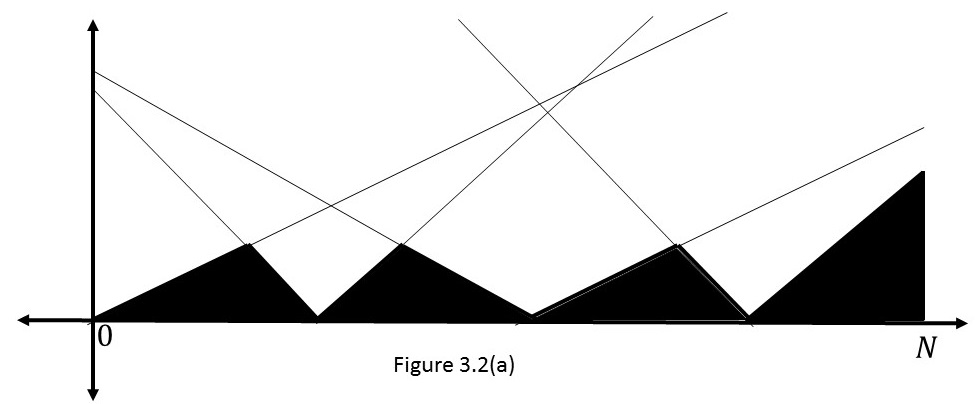}\\
\includegraphics[height = 3.7cm, width = 14cm]{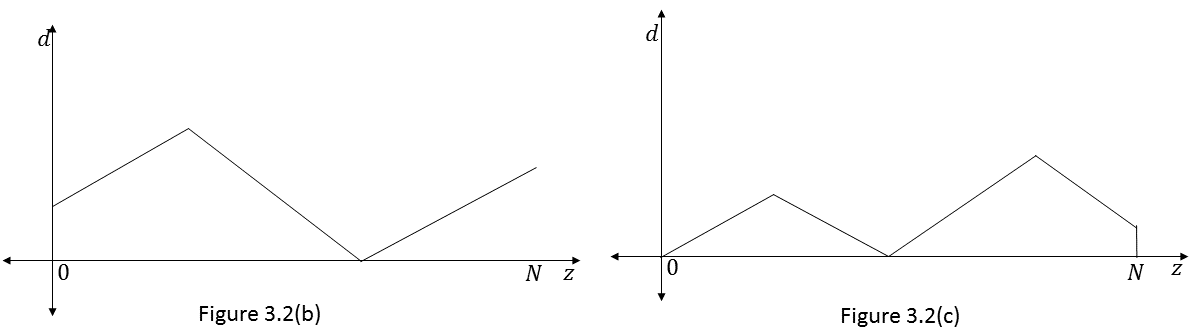}\\
\includegraphics[height = 5.4cm, width = 14cm]{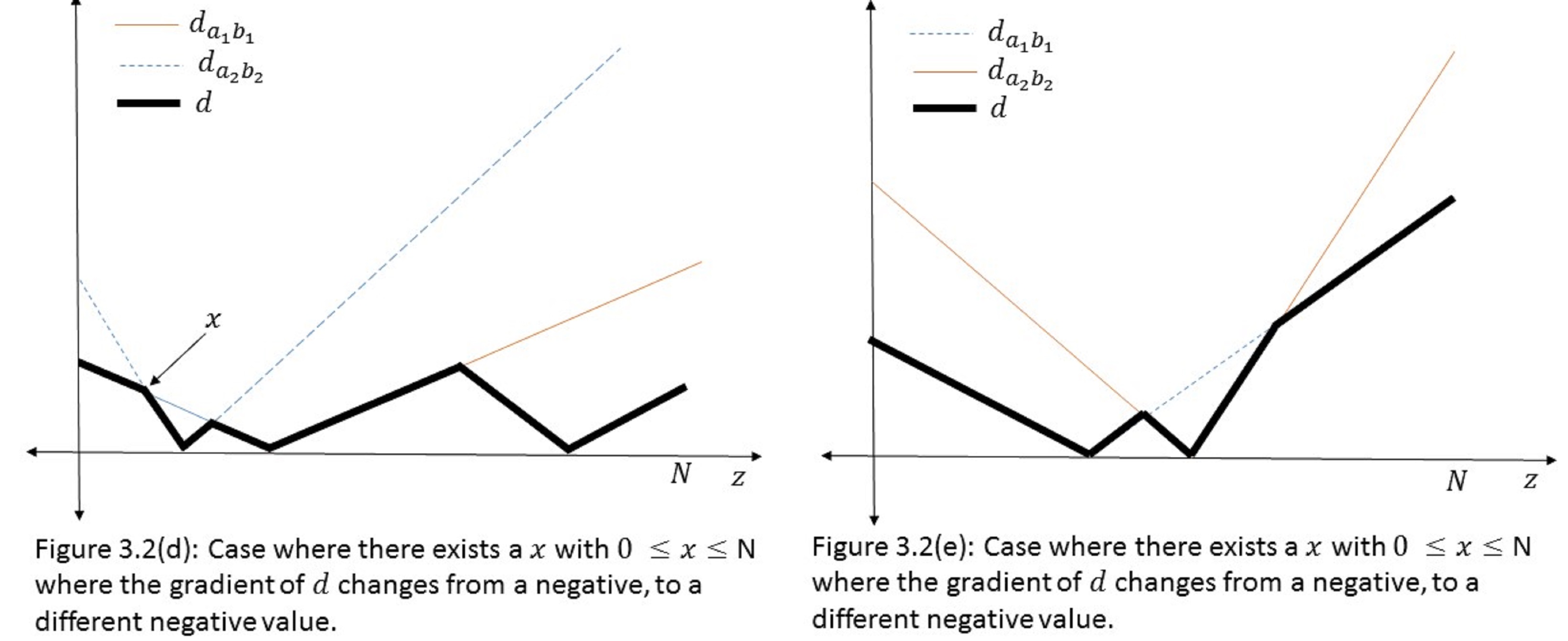} \\

First consider the case where $S_{1}$ doesn't hold. Then $d(0) > 0$, and $d\ensuremath{'}(0) > 0$. Let $d_{C0}$ be the function which is smallest out out of all\\ $\left\{d_{1,1}(z),d_{1,2}(z),..., d_{1,p_{1}}(z),d_{2,1}(z),...,d_{2,p_{2}}(z)\right\}$ at $z = 0$. Hence $d_{C0} = |C + gz|$, where $C > 0$, and $g$ is $1$ or $2$, and $d(z) = d_{C}(z)$ for all $z$ between $0$ and $p$, for some $p < N$. If we change $d_{C0}$ so that $d_{C0} = |gz|$, and define all other functions apart from $d_{C0}$ in the same way, then the area under $d$ decreases by $Cp$ over the region $[0,p]$ (see figure 3.2(f)), while the area under $d$ over $[p,N)$ will be no larger than it was before $d_{C0}(z)$ was redefined to be $|gz|$. Hence the area under $d$ decreases by at least $Cp$ over $[0,N)$. Therefore, in every case where $S_{1}$ doesn't hold, we can decrease the area under $d$ by redefining $d$ so that $S_{1}$ holds \textbf{(}\textsc{Prop1}\textbf{)}.\\ \includegraphics[height = 6cm, width = 14cm]{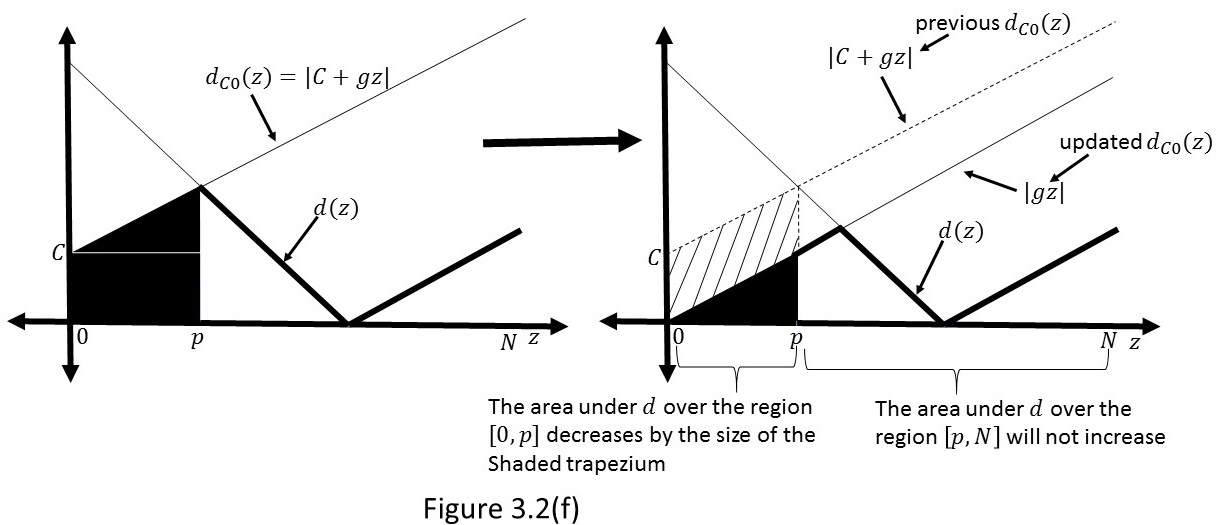} \\
\includegraphics[height = 5.5cm, width = 14cm]{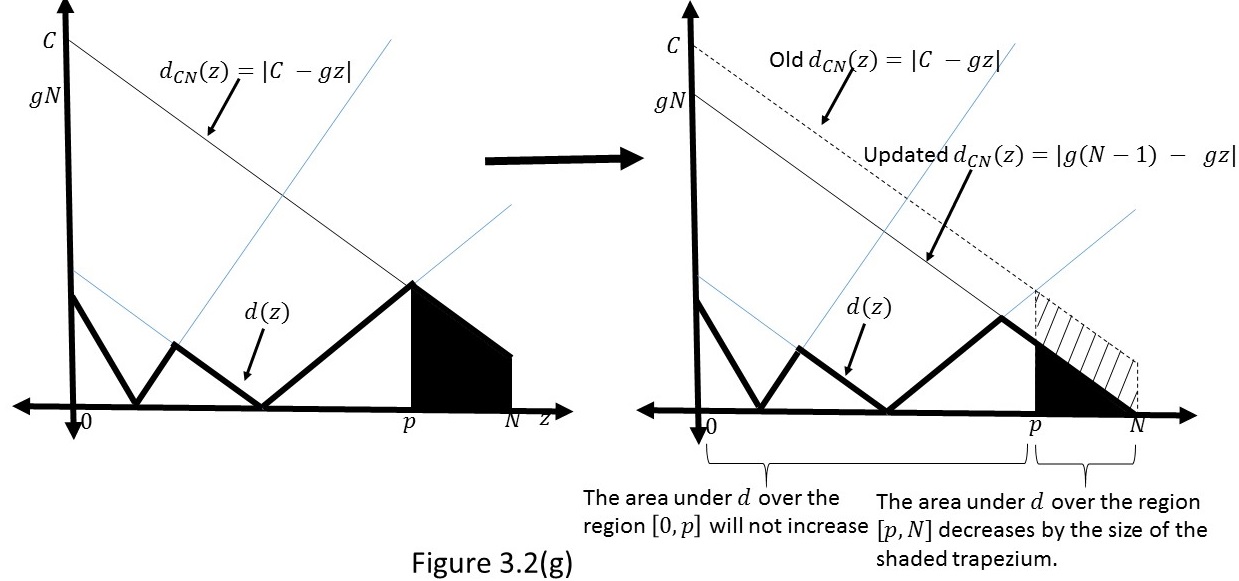} \\
Similarly, in the case where $S_{2}$ doesn't hold, $d(N-1) = C > 0$ and $d\ensuremath{'}(N-1) < 0$. If $d_{CN}$ is the function such that $d_{CN}(N-1) = d(N-1)$, then $d_{CN}(z) = |C - gz|$, where $C > g(N-1)$ (since $d_{CN}(N-1) > 0$). Hence if we redefine $d_{CN}$ such that $d_{CN}(z) = |g(N-1) - gz|$, if $p$ is defined such that $d_{CN}(z) = d(z) \ \forall \ \ p \leq z < N$ (before $d_{CN}(z)$ was defined to be $|g(N-1) - gz|$), then the area under $d$ decreases by at least $(N-1-p)(C-gN)$ (see figure 3.2(g)). Therefore, in every case where $S_{2}$ doesn't hold, we can decrease the area under $d$ by redefining $d$ so that $S_{2}$ holds \textbf{(}\textsc{Prop2}\textbf{)}.\\% \includegraphics[height = 4.5cm, width = 13cm]{fig32G.jpg} \\
\includegraphics[height = 5.5cm, width = 14cm]{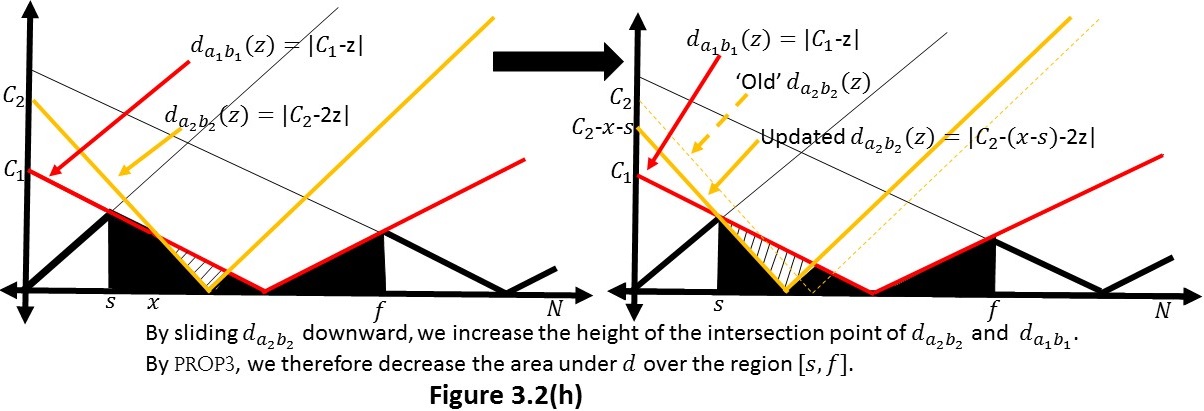} \\
Now consider the case where $S_{3}$ doesn't hold. Then there exists an $x$ with $0 \leq x < N$, and functions $d_{a_{1},b_{1}}$ and $d_{a_{2},b_{2}}$ such that the function which is smallest (so the function which $d$ equals) changes from  $d_{a_{1},b_{1}}$ to $d_{a_{2},b_{2}}$ at $x$, and $d_{a_{1},b_{1}}\ensuremath{'}(x)$ and $d_{a_{2},b_{2}}\ensuremath{'}(x)$ have the same sign.\\ In the case where  $d_{a_{1},b_{1}}\ensuremath{'}(x) < 0$ and $d_{a_{2},b_{2}}\ensuremath{'}(x) < 0$, (see figure 3.2(h)) $d_{a_{2},b_{2}}\ensuremath{'}(x) < d_{a_{1},b_{1}}\ensuremath{'}(x)$ (since $d_{a_{2},b_{2}}(z) > d_{a_{1},b_{1}}(z)$ for $z < x$ with $z \approx x$, and $d_{a_{2},b_{2}}(z) < d_{a_{1},b_{1}}(z)$ for $z > x$ with $z \approx x$. Hence $d_{a_{2},b_{2}}(x)\ensuremath{'} = -2$ and  $d_{a_{1},b_{1}}(x)\ensuremath{'} = -1$. Therefore $d_{a_{2},b_{2}}(z) = |C_{2} - 2z|$, and $d_{a_{1},b_{1}}(z) = |C_{1} - z|$ for some $C_{1},C_{2} > 0$. Now I will let $[s,f]$ be the region such that for all $z$ where $s \leq z \leq f$, $d(z) = d_{a_{1},b_{1}}(z)$ or $d(z) = d_{a_{2},b_{2}}(z)$, and $h$ be such that $h = d_{a_{1},b_{1}}(x) = d_{a_{2},b_{2}}(x)$. Now if we fix $d_{a_{1},b_{1}}$ (and all other functions in $\left\{d_{1,1}(z),d_{1,2}(z),...,d_{1,p_{1}}(z),d_{2,1}(z),...,d_{2,p_{2}}(z)\right\}$), and allow $d_{a_{2},b_{2}}$ to vary (by varying the constant $C_{2}$), then assuming $d_{a_{1},b_{1}}$ and $d_{a_{2},b_{2}}$ intersect somewhere on the region $[s,f]$ when both $d_{a_{1},b_{1}}\ensuremath{'} < 0$ and $d_{a_{2},b_{2}}\ensuremath{'} < 0$, then the area under $d$ over the region $[s,f]$ is determined purely by the height of this intersection point (Mathematically, it can be easily shown that if we define $H$ to be the height of the intersection point of $d_{a_{1},b_{1}}$ and $d_{a_{2},b_{2}}$ when $d_{a_{1},b_{1}}$ and $d_{a_{2},b_{2}}$ are decreasing, and $A1_{s,f}$ to be the area under $d_{a_{1},b_{1}}$ over $[s,f]$, then the area under $d$ over $[s,f]$ is $A1_{s,f}-2H^{2}/9$). Now suppose we redefine $d_{a_{2},b_{2}}$ so that $d_{a_{2},b_{2}}(z) = |C_{2}-(x-s)-2z|$ (as in figure 3.2(h)). Then $d_{a_{1},b_{1}}$ and $d_{a_{2},b_{2}}$ intersect with negative gradient when $z=s$, and the height of this intersection point is $h+(x-s)$.
The following basic lemmas will be very useful for the remainder this proof. 
\begin{lemma3231}
 If $s>0$, $d\ensuremath{'}(z) > 0$ for $z < s$, with $z \approx s$. \end{lemma3231}
\begin{proof}
 At any $z$, $d\ensuremath{'}(z)$ is either $-2,-1,1$ or $2$. If $d\ensuremath{'}(z) < 0$ for $z < s$ with $z \approx s$, then since $d_{a_{1},b_{1}}(z) = d(z)$ for $z > s$ with $z \approx s$, $d_{a_{1},b_{1}}$ would still be the smallest function for $z < s$, with $z \approx s$. But this is a contradiction to $[s,f]$ being the region which either $d_{a_{1},b_{1}}$ or $d_{a_{1},b_{1}}$ are the smallest functions over. \end{proof}
\begin{lemma3232} In the case where the intersection point of $d_{a_{1},b_{1}}$ and $d_{a_{2},b_{2}}$ when $d_{a_{1},b_{1}}\ensuremath{'},d_{a_{2},b_{2}}\ensuremath{'} < 0$ is when $z = s$, $d$ satisfies $S_{3}$ over $[s,f]$. \end{lemma3232}
\begin{proof} After $d_{a_{2},b_{2}}$ is redefined, $d_{a_{1},b_{1}}$ and $d_{a_{2},b_{2}}$ are clearly still the smallest functions over $[s,f]$. Hence we only need to consider the intersection points of $d_{a_{1},b_{1}}$ and $d_{a_{2},b_{2}}$ to prove the lemma. By the nature of the functions $d_{a_{1},b_{1}}$ and $d_{a_{2},b_{2}}$, $d_{a_{1},b_{1}}$ and $d_{a_{2},b_{2}}$ can have at most one intersection point when both $d_{a_{1},b_{1}}\ensuremath{'}$ and $d_{a_{2},b_{2}}\ensuremath{'}$ have the same sign. Now $d_{a_{1},b_{1}}$ and $d_{a_{2},b_{2}}$ intersect when both have negative gradient when $z = s$. Hence the only $z$ where $d$'s gradient changes, but the sign of $d$'s gradient remains the same over $[s,f]$, is when $z = s$. But Lemma 3.2.3.1 implies that the gradient of $d$ changes from a positive to a negative value when $z = s$. Hence $S_{3}$ is satisfied over $[s,f]$. \end{proof} 

As a result, we can conclude $x > s$, since otherwise $S_{3}$ would hold over $[s,f]$ when $d_{a_{2},b_{2}}(z)$ was $|C_{2}-2z|$. Hence the height of the intersection point of $d_{a_{1},b_{1}}$ and $d_{a_{2},b_{2}}$ when both $d_{a_{1},b_{1}}\ensuremath{'},d_{a_{2},b_{2}}\ensuremath{'} < 0$ is increased when $d_{a_{2},b_{2}}$ is redefined, so the area under $d$ over $[s,f]$ is decreased by redefining $d_{a_{2},b_{2}}$ in such a way that $d$ satisfies $S_{3}$ over $[s,f]$. Since the value of $d(z)$ doesn't increase when $d_{a_{2},b_{2}}$ is redefined for $0 \leq z < s$ and $f < z < N$, the area under $d$ over all other regions apart from $[s,f]$ does not increase when $d_{a_{2},b_{2}}$ is redefined. Hence the area under $d$ over $[0,N)$ is decreased by redefining $d$ so that $S_{3}$ is satisfied over $[s,f]$.\\
A very similar argument can show that in the case where there exists a $z$ such that the function which is smallest changes from  $d_{A_{1},B_{1}}$ to $d_{A_{2},B_{2}}$ at $z$, and $d_{A_{1},B_{1}}\ensuremath{'},d_{A_{2},B_{2}}\ensuremath{'} > 0$, then we can redefine $d_{A_{2},B_{2}}$ so that the area under $d$ over $[0,N)$ is decreased and $S_{3}$ is satisfied over $[S,F]$ (where $S$ and $F$ are defined such that $d$ equals either $d_{A_{1},B_{1}}$ or $d_{A_{2},B_{2}}$ for all $z$ with $S \leq z \leq F$).\\ Therefore, in every case where there exists $z$ such that the gradient of $d$ changes but keeps the same sign at $z$, we can decrease the area under $d$ over $[0,N)$ by redefining one of the functions in $\left\{d_{1,1}(z),d_{1,2}(z),...,d_{1,p_{1}}(z),d_{2,1}(z),...,d_{2,p_{2}}(z)\right\}$, so that $d$ satisfies $S_{3}$ over each of the intervals $[s,f]$ and $[S,F]$. Hence, in such a case we can decrease the area under $d$ while making $d$ satisfy $S_{3}$ over $[0,N)$ \textbf{(}\textsc{Prop3}\textbf{)}.\\ By combining \textsc{Prop1},\textsc{Prop2}, and \textsc{Prop3}, we can conclude that in every case where $d$ is not a sequence of triangles (so at least one of $S_{1}$,$S_{2}$ or $S_{3}$ doesn't hold for $d$), we can decrease the area under $d$ by redefining $d$ so that $d$ is a sequence of triangles (so $S_{1}$,$S_{2}$ and $S_{3}$ all hold for $d$ on $[0,N)$). \end{proof}

In any five kangaroo method, we are given a set of functions of the same form as those in lemma 3.2.3. To minimise the expected number of group operations until the closest useful pair collides, we need to minimise the area under $d$, by choosing the $C_{i}$ and $C_{j}$s appropriately in functions of the form of lemma 3.2.3. Hence we may assume that when the area under $d$ is as small as possible, $d$ will be a sequence of triangles. Therefore, since this theorem is stating a lower bound on the expected number of group operations until the closest useful pair collides, for the purposes of the proof of this theorem, $d$ can be assumed to be a sequence of triangles.\\The following lemma will be useful in finding how small the area under $d$ can be, given that $d$ can be assumed to be a sequence of triangles.

\begin{lemma324} Fix $n,N \in \mathbb{N}$, and let the gradients $G_{1}$,$G_{2}$,...,$G_{n} \in \mathbb{R}\setminus{\left\{0\right\}}$ be fixed. Let $R_{1}$,$R_{2}$,...,$R_{n}$ be such that $R_{i} \geq 0$, $\sum_{i=1}^{n} R_{i} = N$, and such that the sum of areas of triangles of base $R_{i}$ and gradient $G_{i}$ is minimised. Then all triangles have the same height.\end{lemma324}

\begin{proof}[Proof of Lemma 3.2.4] I will label the triangles such that the $i^{th}$ triangle ($T_{i}$) is the triangle which has $i-1$ triangles to the left of it. $G_{i}$ can be considered to be the gradient of the slope of $T_{i}$, and $R_{i}$ can be considered to be the size of the region which $T_{i}$ occupies. Also let $A_{T}$ be the sum of area under these triangles. A diagram of the situation is shown in Figure 3.1.

\begin{figure}
\includegraphics[width=12cm,height=5cm,left]{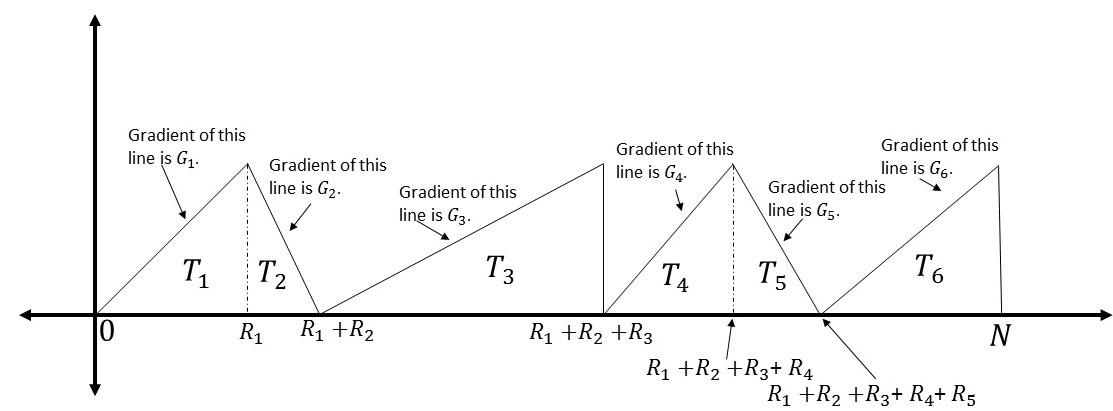}
\caption{Diagram of the kind of situation being considered in lemma 3.2.4}
\end{figure}

Then we have $$\sum_{i=1}^n R_{i} = N$$ and $$A_{T} = \sum_{i=1}^n \frac{|G_{i}|.R_{i}^2}{2}$$ In the arrangement where the heights of all the triangles are the same, for each $1 \leq i < j < n$, $R_{i}/R_{j} = |G_{j}|/|G_{i}|$. Suppose we change the ratio of $R_{i}$ to $R_{j}$, so that the new $R_{i}$ is $R_{i} + \epsilon$, and $R_{j}$ becomes $R_{j} - \epsilon$ (note that $\epsilon$ can be greater than or less than $0$). Then the area of $T_{i}$ becomes $(G_{i}(R_{i} + \epsilon)^2)/2 = (G_{i}(R_{i}^2 + 2R_{i}\epsilon + \epsilon^2))/2$ while the area of $T_{j}$ becomes $G_{j}(R_{j}^2 - 2R_{j}\epsilon + \epsilon^2)/2$. $A_{T}$ therefore changes by \begin{center}

$$\frac{|G_{i}|(R_{i}+\epsilon)^2}{2} + \frac{|G_{j}|(R_{j}-\epsilon)^2}{2} - \frac{|G_{i}|R_{i}^2}{2} - \frac{|G_{j}|R_{j}^2}{2}$$ $$= \frac{|G_{i}|R_{i}^2}{2} + \frac{|G_{j}|R_{j}^2}{2} + |G_{i}|R_{i}\epsilon - |G_{j}|R_{j}\epsilon - \frac{|G_{i}|R_{i}^2}{2} - \frac{|G_{j}|R_{j}^2}{2} + \frac{\epsilon^2}{2} + \frac{\epsilon^2}{2}$$ $$= |G_{i}|R_{i}\epsilon - |G_{j}|R_{j}\epsilon + \epsilon^2 = \epsilon^2 > 0$$
\end{center}

Therefore, adjusting the ratio of the size of any two regions away from that which ensures the height of the triangles is the same, always strictly increases the sum of the area of the triangles. It follows that when the heights of all the triangles are equal, the sum of the area beneath these triangles is minimised. \end{proof}

\begin{lemma325} In the case where $d$ is a sequence of triangles, each useful pair is either never the closest useful pair, or it is the closest useful pair over a single region. \end{lemma325}

\begin{proof} By a useful pair ($P$) being the closest useful pair only over a single region, I mean that there is a single interval $[s,f]$, with $0 \leq s < f < N$, such that $P$ is the closest useful pair for all $z$ where $s \leq z \leq f$. This is in contrast to there being intervals $[s_{1},f_{1}]$, and $[s_{2},f_{2}]$, where $0 \leq s_{1} < f_{1} < N$, and $f_{1} < s_{2} < f_{2} < N$, such that $P$ is the closest useful pair for all $z$ with $s_{1} \leq z \leq f_{1}$, and $s_{2} \leq z \leq f_{2}$, but $P$ is not the closest useful pair for all $z$ where $f_{1} < z < s_{2}$.\\ 
The result of this lemma follows easily from the fact that if $d$ is a sequence of triangles, then for every $z$ where the closest useful pair changes (say from $P_{1}$ to $P_{2}$), the gradients of the initial distance functions between $P_{1}$, and $P_{2}$ have opposite sign. \end{proof}

\begin{proof}[Proof of Theorem 3.2]
The lemmas can now be used to prove the theorem. Firstly, label all useful pairs in any way from $1$ to $n$, where $n$ is the number of useful pairs of kangaroos in $A$. Now since $d$ can be assumed to be a a sequence of triangles, lemma 3.2.5 implies that a useful pair is either never the closest useful pair, or there exists a single region for which a useful pair is closest over. For each $i$, in the case where pair $i$ is the closest useful pair over some region, let $R_{i}$ denote the single region which pair $i$ is the closest useful pair over, and $R_{i_{s}}$ denote the size of $R_{i}$. In the case where pair $i$ is never the closest useful pair, define $R_{i_{s}}$ to be $0$. Then $\sum_{i=1}^n R_{i_{s}} = N$. Letting $d_{i}$ denote the distance function between pair $i$ for each $i$, by Lemma 3.2.2, $d_{i}(z)$ is of the form $|C \pm g_{i}z|$, where $g_{i}$ is 1 or 2, and $C \in \mathbb{R}$. It is clear from the diagram in section 2.4 that such functions feature at most two triangles, both of which have slopes of gradient with an absolute value of $g_{i}$. Hence $d_{i}$ must feature at most two such triangles over $R_{i}$. Now for every $i$, pair $i$ is the closest useful pair over $R_{i}$, so $d(z) = d_{i}(z) \ \forall z \in R_{i}$. Therefore, $d$ features at most two triangles over $R_{i}$, both of which have slopes that have a gradient with an absolute value of $g_{i}$. Therefore, if we let $T_{g_{1}}$ and $T_{g_{2}}$ denote the number of triangles in $d$ where the absolute values of their gradients are 1 and 2 respectively, $T_{g_{1}} \leq 2p_{g_{1}} = 2N_{T}(N_{W1} + N_{W2})$ and $T_{g_{2}} \leq 2p_{g_{2}} = 2N_{W1}N_{W2}$. Now by Lemma 3.2.4, for the  area under $d$ to be minimised, all triangles must have the same height. This implies that all triangles with the same gradient in $d$ must cover a region of the same size. Defining $R_{T_{1}}$ and $R_{T_{2}}$ to denote the size of the regions covered by each triangle of gradients 1 and 2 respectively, Lemma 3.2.4 also implies that $R_{T_{1}} = 2R_{T_{2}}$. Then since all triangles cover the domain of $d$, $N = T_{g_{1}}R_{T_{1}} + T_{g_{2}}R_{T_{2}}$ $\leq 2N_{T}(N_{W1} + N_{W2})R_{T_{1}} + 2N_{W1}N_{W2}R_{T_{2}}$ $= 4N_{T}(N_{W1} + N_{W2})R_{T_{2}} + 2N_{W1}N_{W2}R_{T_{2}}$. Hence $R_{T_{2}} \geq N/(4N_{T}(N_{W1} + N_{W2}) + 2N_{W1}N_{W2})$. Now since all triangles in $d$ have the same height, the average of $d$ over all $z$ with $0 \leq z < N$ is half the height of all triangles. Therefore, since the height of a triangle of gradient 2 is $2R_{T_{2}}$, $R_{T_{2}}$ gives the average of $d$, which is the average distance between the closest useful pair. Then by Lemma 3.2.1, the expected number of group operations until the closest useful pair collides in the case where the area under $d$ is minimised (and hence the average distance between the closest useful pair is minimised) is $10\sqrt{R_{T_{2}}} \geq 10\sqrt{N/(4N_{T}(N_{W1} + N_{W2}) + 2N_{W1}N_{W2})}$. 
\end{proof}

By plugging in various values of $N_{T}$, $N_{W1}$ and $N_{W2}$ into this formula under the constraint that $N_{T} + N_{W1} + N_{W2} = 5$, we can gather a lower bound on the expected number of group operations until the closest useful pair collides, when different numbers of each type of walk are used. The following table shows the methods which achieved the three best lower bounds. In the table, a tuple of the form $(x,y,z)$ denotes an algorithm that uses $x$ \textsc{tame} kangaroos, $y$ \textsc{wild1} kangaroos, and $z$ \textsc{wild2} kangaroos.

\begin{tabular}{|l|p{6cm}|}
\hline
$10\sqrt{\frac{N}{(4N_{T}(N_{W1}+N_{W2}) + 2N_{W1}N_{W2})}}$ & Number of walks of each type used \\ \hline
$1.8898\sqrt{N}$ & $(2,1,2),(2,2,1)$ \\ \hline
$1.9611\sqrt{N}$ & $(3,1,1)$ \\ \hline
$2.0412\sqrt{N}$ & $(1,2,2),(2,0,3),(2,3,0)$,\\ & $(3,0,2),(3,2,0)$ \\ \hline

\end{tabular} 

From this, it can be seen that a $(2,2,1)$ 5 kangaroo method has the minimal lower bound on the expected number of group operations until the closest useful pair collides, of $1.8898\sqrt{N}$ group operations. A $(2,1,2)$ method doesn't need to be considered as a separate case, since this is clearly equivalent to a $(2,2,1)$ method. If one starts two \textsc{wild1} kangaroos at $h$ and $g^{0.7124N}h$, a \textsc{wild2} kangaroo at $g^{1.3562N}h^{-1}$, and two \textsc{tame} kangaroos at $g^{0.9274N}$ and $g^{0.785N}$, the lower bound $(2,2,1)$ method of $1.8898\sqrt{N}$ group operations is realised. % (The approach that was used to find these starting positions is very similar to the one I used to the find the starting positions that minimise the number of group operations required to solve the IDLP. The pseudocode for this algorithm can be found on which is outlined in Algorithm 2 on page...) I used to find this lower bound in [section on how I found optimal overall layout]
The table shows that in any other 5 kangaroo method, the closest useful pair can't collide in less than $1.9611\sqrt{N}$ group operations on average.

Since the closest useful pair of kangaroos is by far the most significant in determining the running time of any kangaroo algorithm (for instance, in the 3 kangaroo method the expected number of group operations until the closest useful pair collides could be as low as $1.8972\sqrt{N}$ group operations, while the expected number of group operations until any pair collides can only be as low as $1.818\sqrt{N}$ group operations), a method that uses 2 \textsc{tame}, 2 \textsc{wild1}, and 1 \textsc{wild2} kangaroos is most likely to be the optimal 5 kangaroo method, out of all methods that only use \textsc{tame},\textsc{wild1}, and \textsc{wild2} kangaroos.

\section[Designing a $(2,2,1)$ kangaroo algorithm]{Where the Kangaroos should start their walks, and what average step size should be used?}

Given that, in any 5 kangaroo method that uses only \textsc{tame}, \textsc{wild1}, and \textsc{wild2} kangaroos one should use 2 \textsc{wild1}, 1 \textsc{wild2}, and 2 \textsc{tame} kangaroos, the next question to consider is where abouts the kangaroos should start their walks. In the analysis that answers this question, the question of what average step size to use will be answered also. I will now state some definitions and remarks that will be used throughout the remainder of this thesis.\\ 

\begin{itemize} \item \textbf{Remark 1: } Firstly, I will redefine the IDLP to be to find $z$, given $h = g^{zN}$, when we're given $g$ and $h$, and that $0 \leq z < 1$.
\item \textbf{Remark 2: } At this stage, I will place one constraint on the starting positions of all kangaroos. This being, that all \textsc{Wild1},\textsc{Wild2}, and \textsc{Tame} kangaroos will respectively start their walks at positions of the form $aN + zN$, $bN - zN$, and $cN$, for universal constants $a$,$b$, and $c$, that are independent of the interval size $N$. This is the only constraint I will place on the starting positions at this stage.
\item \textbf{Remark 3: } I will let $D_{i,z}$ denote the initial distance between the $i^{th}$ closest useful pair of kangaroos for some specified $z$. It follows from Remark 2 that $D_{i,z} = d_{i,z}N$, for some $d_{i,z}$ independent of $N$. 
\item \textbf{Remark 4: } The average step size $m$ will be defined to be $c_{m}\sqrt{N}$, for some $c_{m}$ independent of $N$.
\item \textbf{Remark 5: } $S_{i,z}$ will be defined such that $S_{i,z}$ denotes the expected number of steps the back kangaroo requires to catch up to the front kangaroo's starting position in the $i^{th}$ closest useful pair, for some specified $z$. It is clear that $S_{i,z} = \lceil \frac{D_{i,z}}{m}\rceil$. Since $D_{i,z} = d_{i,z}N$, and $m = c_{m}\sqrt{N}$ for some $d_{i,z}$, and $c_{m}$ which are independent of $N$, we can say $S_{i,z} = \lceil s_{i,z}\sqrt{N} \rceil$ for some $s_{i,z}$ independent of $N$. For the typical interval sizes over which one uses kangaroo methods to solve the IDLP ($N > 2^{30}$), this can be considered to be $s_{i,z}\sqrt{N}$.
\item \textbf{Remark 6: } $C_{i,z}$ and $c_{i,z}$ will be defined such that $C_{i,z} = \sum_{j=1}^{i} S_{i,z}$, and $c_{i,z}\sqrt{N} = C_{i,z}$. One can see from the definition of $S_{i,z}$ that $c_{i,z}$ is independent of the interval size $N$. \end{itemize}

I will answer the question that titles this section by first presenting a formula that can compute the running time of any (2,2,1) 5-kangaroo algorithm (see section 3.3.1), and then by showing how this formula can be used to find the best starting positions and average step size to use (see section 3.3.2). 

\newpage

\subsection{Formula for computing the running time of a (2,2,1)-5 kangaroo algorithm}
\begin{theorem331} In any 5 kangaroo method which uses 2 \textsc{tame}, 1 \textsc{wild2}, and 2 \textsc{wild1} kangaroos, the expected number of group operations required to solve the IDLP is approximately \\ $5\left(\left(\int_{0}^1 c_{z} dz\right) + o(1)\right)\sqrt{N} + O(\log(N))$, where \[c_{z} = \sum\limits_{i=1}^{8}\left(e^{\frac{(-is_{i,z} + c_{i,z})}{c_{m}}}(\frac{c_{m}}{i} + s_{i,z}) -  e^{\frac{-is_{i+1,z} + c_{i,z}}{c_{m}}}(\frac{c_{m}}{i} + s_{i+1,z})\right)\]. \end{theorem331} %THIS FORMULA WILL NEED TO BE CHANGED TO ACCOMMODATE FOR THE FACT THAT THE DEFN OF S_{I,Z} HAS CHANGED (NEW SIZ - OLD SIZ + 1
\begin{proof}
Any 5 kangaroo method can be broken into the following 4 disjoint stages;
\begin{itemize}
\item \textbf{Stage 1-} The stage where the algorithm is initialised. This involves computing the starting positions of the kangaroos, assigning a step size to each group element, and computing and storing the group elements which kangaroos are multiplied by at each step (so computing $g^{H(x)}$ for each $x$ in $G$).
\item \textbf{Stage 2-} The period between when the kangaroos start their walks, and the first collision between a useful pair occurs.
\item \textbf{Stage 3-} The period between when the first collision occurs, and when both kangaroos have visited the same distinguished point.
\item \textbf{Stage 4-} The stage where $z$ is computed, using the information gained from a useful collision.
\end{itemize}

\subsubsection{Number of group operations required in Stage 1}

To find the starting positions of the kangaroos, we require one inversion (to find the starting position of the kangaroo of type \textsc{wild2}), 3 multiplications (in finding the starting positions of all wild kangaroos), and 5 exponentiations (one to find the starting position of each kangaroo). All of these operations are $O(\log(N))$ in any group. As explained in [8], the step sizes can be assigned in $O(\log(N))$ time also using a hash function. To pre-compute the group elements which kangaroos can be multiplied by at each step, we need to compute $g^{s}$, for each step size $s$ that can be assigned to a group element (this is the number of values which the function $H$ can take). The number of step sizes used in kangaroo methods is generally between 20 and 100. This was suggested by Pollard in [8]. Applying less than a constant number (100) of exponentiation operations requires a constant number of group operations (so independent of the interval size).\\ Summing together the number of group operations required by each part of Stage 1, we see that the number of group operations required in Stage 1 is $O(\log(N))$.

\subsubsection{Number of group operations required in Stage 2}

%To analyse the number of group operations required in stage 2, I will make the assumption that the expected number of group operations until the 
To analyse the number of group operations required in stage 2, I will define $Z(z)$ to be a random variable on $\mathbb{N}$ such that $\Pr \big(Z(z) = k\big)$ denotes the probability that for some specified $z$, the first collision occurs after $k$ steps. From this, one can see that $\mathrm{E}\big(Z(z)\big) = \sum_{k=0}^{\infty} k\Pr\big(Z(z) = k\big)$, gives the expected number of steps until the first collision occurs, at our specified $z$.

\subsubsection*{Important Remark} To compute $\mathrm{E}(Z(z))$, I will compute the expected number of steps until the first useful collision occurs in the case where in every useful pair of kangaroos, the back kangaroo takes the expected number of steps to catch up to the front kangaroos starting position, and make the assumption that this is proportional to the expected number of steps until the first useful collision occurs across all possible random walks. This assumption was used implicitly in computing the running time of the three kangaroo method in [3], and is a necessary assumption to make, since calculating the expected number of steps until the first useful collision occurs across all possible walks is extremely difficult.\\ The following lemma will be useful in computing $\mathrm{E}(Z(z))$. 

\begin{lemma332} In the case where in every useful pair of kangaroos, the back kangaroo takes the expected number of steps to catch up to the starting position of the front kangaroo, for every $k$, $\Pr (Z(z) = k) = \sum\limits_{j=1}^i {i \choose j} \frac{1}{m^j}e^{\frac{-ik - i + C_{i,z} + j}{m}}$ where $i$ is the number of pairs of kangaroos for which the back kangaroo has caught up to the front after $k$ steps. \end{lemma332}

\begin{proof}[Proof] Let $k \in \mathbb{N}$, and $i$ be such that $S_{i}(z) \leq k < S_{i+1}(z)$. In the case where in every pair of kangaroos, the back kangaroo takes the expected number of steps to catch up to the starting position of the front kangaroo, when $S_{i}(z) \leq k < S_{i+1}(z)$, there will be exactly $i$ pairs where the back kangaroo will be walking over a region that has been traversed by the front kangaroo (these will be the $i$ closest useful pairs). Therefore, exactly $i$ pairs can collide after $k$ steps, for all $S_{i}(z) \leq k < S_{i+1}(z)$. In order for the first collision to occur after exactly $k$ steps, we require that the $i$ pairs that can collide avoid each other for the first $k-1$ steps, and then on the $k^{th}$ step, $j$ pairs collide for some $j$ between $1$ and $i$. I will define $E_{k,j}$ to be the event that there are no collisions in the first $k-1$ steps, and then on the $k^{th}$ step, exactly $j$ pairs collide. Since $E_{k,j}$ and $E_{k,l}$ are disjoint events for $j \neq l$, we can conclude that $\Pr\big(Z(z) = k\big) = \sum\limits_{j=1}^i P(E_{k,j})$ \textbf{(1)}. Now $\Pr (E_{k,j})$ can be computed as follows;
In any instance where $E_{k,j}$ occurs, we can define sets $X$ and $Y$ such that $X$ is the set of all $x$ where the $x^{th}$ closest useful pair of kangaroos doesn't collide in the first $k-1)$ steps, but does collide on the $k^{th}$ step, and $Y$ is the set of all $y$ where the $y^{th}$ closest useful pair doesn't collide in the first $k$ steps. %$X = \left\{ {x \in {1,2,…,i}|the x^{th} closest pair collides after k steps} \right\{$, and $Y = \left\{ {y \in {1,2,….,i}|the y^{th} closest pair doesn’t collide in the first y steps}\right\}$.
Now in Stage 2 of the van Oorshot and Weiner method, I explained how at any step, the probability that a pair collides once the back kangaroo has caught up to the path of the front is $1/m$. Hence for any $x \in X$, the probability that the back kangaroo in the $x^{th}$ closest useful pair avoids the path of the front kangaroo for the first $k-1$ steps, but lands on an element in the front kangaroos walk on the $k^{th}$ step is $\frac{1}{m}(1 - \frac{1}{m})^{k-S_{x,z}}\approx \frac{1}{m}e^{\frac{-k+S_{x,z}}{m}}$, while for any $y \in Y$, the probability that the $y^{th}$ closest useful pair doesn’t collide in the first $k$ steps is $(1-\frac{1}{m})^{k+1-S_{y,z}} \approx e^{\frac{-k-1+S_{y,z}}{m}}$. Now before any collisions have taken place, the walks of any 2 pairs of kangaroos are independent of each other. Therefore, the probability that the pairs in $X$ all first collide on the $k^{th}$ step, while all the pairs in $Y$ don't collide in the first $k$ steps is $\prod_{x \in X} \frac{1}{m}e^{\frac{-k+S_{x,z}}{m}} \prod_{y \in Y} e^{\frac{-k-1+S_{y,z}}{m}}$ $= \frac{1}{m^j}e^{\frac{-jk + \sum_{x \in X} S_{x,z}}{m}}e^{\frac{-(i-j)k - (i-j) + \sum_{y \in Y} S_{y,z}}{m}}$ $= \frac{1}{m^j}e^{\frac{-ik -i + j + C_{i,z}}{m}}$. Now since there are ${i \choose j}$ ways for $j$ out of the $i$ possible pairs to collide on the $k^{th}$ step, we obtain the formula $P(E_{k,j}) = {i \choose j} \frac{1}{m^j}e^{\frac{-ik -i + j + C_{i,z}}{m}}$. When substituting this result back into \textbf{(1)}, we obtain the required result of $\Pr (Z(z) = k) = \sum_{j=1}^{i} {i \choose j} \frac{1}{m^j}e^{\frac{-ik -i + j + C_{i,z}}{m}}$\end{proof} 
I will now define $p_{z}$ to be the function such that $p_{z}(k) = k\Pr (Z(z) = k)$, for all $k \in N$. Hence $\mathrm{E}(Z(z)) = \sum_{k=0}^{\infty} p_{z}(k)$. Now in a 5 kangaroo method that uses 2 \textsc{Tame}, 2 \textsc{Wild1}, and 1 \textsc{Wild2} walks, there there are 8 pairs that can collide to yield a useful collision (4 \textsc{Tame}/\textsc{Wild1} pairs, 2 \textsc{Tame}/\textsc{Wild2} pairs, and 2 \textsc{Wild1}/\textsc{Wild2} pairs). Hence for any $k$, the number of pairs where the back kangaroo has caught up to the front kangaroos starting position can be anywhere between 0 and 8. %Therefore, $p$ is a piecewise function, defined on the nine intervals $p = k\sum_{j=1}^{i} {i \choose j} \frac{1}{m^j}e^{\frac{-ik -i + j + C_{i,z}}{m}}$ for all $i$ between 0 and 8, and $p$ can be broken into the nine regions [Si+1,Si+1], for i all 0 <= I <= 8 (and where S0 := 0, and S9 := infinity). PLACE THIS IN From now on I will consider p as a continuous function, where the formula that holds for k element of N with k element of [s(i)+1,Si] holds for all k element of R, with k element of [Si+1,Si+1+1). I will define pi to be the function such that pi(k) = p(k) for all k element of [Si+1,Si+1+1).
Therefore, 
\[
  p_{z}(k) = \left\{\def\arraystretch{1.2}%
  \begin{array}{lr}
  0 & \text{$1 \leq k < S_{1,z}$}\\
    \frac{k}{m}e^{\frac{-k+C_{1,z}}{m}} & \text{$S_{1,z} \leq k < S_{2,z}$}\\
    k(\frac{1}{m}e^{\frac{-2k-1+C_{2,z}}{m}} + \frac{1}{m^2}e^{\frac{-2k + C_{2,z}}{m}})& \text{$S_{2,z} \leq k < S_{3,z}$}\\
    \sum\limits_{j=1}^3 k{3 \choose j}\frac{1}{m^j}e^{\frac{-3k - 3 + j + C_{3,z}}{m}} & \text{$S_{3,z} \leq k < S_{4,z}$}\\
    \sum\limits_{j=1}^4 k{4 \choose j}\frac{1}{m^j}e^{\frac{-4k - 4 + j + C_{4,z}}{m}} & \text{$S_{4,z} \leq k < S_{5,z}$}\\
    \sum\limits_{j=1}^5 k{5 \choose j}\frac{1}{m^j}e^{\frac{-5k - 5 + j + C_{5,z}}{m}} & \text{$S_{5,z} \leq k < S_{6,z}$}\\
    \sum\limits_{j=1}^6 k{6 \choose j}\frac{1}{m^j}e^{\frac{-6k - 6 + j + C_{6,z}}{m}} & \text{$S_{6,z} \leq k < S_{7,z}$}\\
    \sum\limits_{j=1}^7 k{7 \choose j}\frac{1}{m^j}e^{\frac{-7k - 7 + j + C_{7,z}}{m}} & \text{$S_{7,z} \leq k < S_{8,z}$}\\
    \sum\limits_{j=1}^8 k{8 \choose j}\frac{1}{m^j}e^{\frac{-8k - 8 + j + C_{8,z}}{m}} & \text{$S_{8,z} \leq k < \infty$}\\
  \end{array}\right.
\]
%
%Following this, see booklet on proof that integral over all k of p(k) – O(1)  <= sum over p(k) <= sum over integral of p(k) + O(1).
%Then see the booklet, ‘end of section on how to analyse the 3 wild 2 tame roo method’. 
For the rest of this thesis, I will consider $p_{z}$ as a continuous function.
The following result will be useful in computing $\mathrm{E}(Z(z))$.
\begin{theorem332} $\int_{1}^{\infty} p_{z}(k)dk - O(1) \leq \mathrm{E}(Z(z)) \leq \int_{1}^{\infty} p_{z}(k) dk + O(1)$.\end{theorem332}

%Place diagram in this proof?
\begin{proof}
The proof of this will use the following lemma.
\begin{lemma3321} $p_{z}(k)$ is $O(1)$ $\forall k$ \end{lemma3321}
\begin{proof}[Proof] Let $0 \leq z < 1$.  Then $\forall \ k < S_{1,z}$, $p_{z}(k) = 0$, while $\forall \ k \geq S_{1,z}$, $p_{z}(k) = k\sum_{j=1}^{i} {i \choose j} \frac{1}{m^j}e^{\frac{-ik -i + j + C_{i,z}}{m}}$, for some $1 \leq i \leq 8$. Hence for the purposes of the proof of this lemma, we can assume $k \geq S_{1,z}$. I will state some facts that will make the argument of this proof flow more smoothly. \

\textbf{Fact 1: } For $k$ such that $S_{i,z} \leq k < S_{i+1,z}$, $e^{\frac{-ik - i + j + C_{i,z}}{m}} \leq 1$.\\ This holds because since $k \geq S_{i,z}$, $ik \geq iS_{i,z} \geq \sum_{j=1}^{i} S_{j,z} = C_{i,z}$. Also, $i \geq j$. Hence $-ik - i + j + C_{i,z} \leq 0$, and $e^{\frac{-ik - i + j + C_{i,z}}{m}} \leq 1$. \

\textbf{Fact 2: } No useful pair of kangaroos can start their walks further than a distance of $6N$ apart on an interval of size $N$. Also, the average step size is at least $0.06697\sqrt{N}$. Both of these facts will be explained in 'Remark regarding Lemma 3.3.2' on Page 44.\ 

\textbf{Fact 3: } The interval size can be assumed to be greater than $2^{30}$. This was explained in the introduction, and can be assumed because one would typically use baby-step giant-step algorithms to solve the IDLP on intervals of size smaller than $2^{30}$.\

I will show that $p_{z}(k)$ is $O(1)$ for two separate cases. \

\textbf{Case 1: } The case where $S_{8,z} \leq m/8$.\\ First, consider the situation in this case where $S_{1,z} \leq k \leq m/8$. Let $i$ be such that $S_{i,z} \leq k < S_{i+1,z}$. Then by Fact 1, Fact 2, Fact 3, and the condition that $k \leq m/8$, $p_{z}(k) = \sum_{j=1}^{i} {i \choose j} \frac{k}{m^j}e^{\frac{-ik -i + j + C_{i,z}}{m}}$ $\leq \sum_{j=1}^{i} {i \choose j} \frac{k}{m^j}$ $\leq \sum_{j=1}^{i} {i \choose j} \frac{\frac{m}{8}}{m^{j}}$ $\leq \sum_{j=1}^{i} {i \choose j} \frac{1}{8(0.06697\sqrt{2^{30})^{j-1}}}$, which is clearly $O(1)$.\\ Now consider the case where $k > m/8$. Then since $S_{8,z} < m/8$, $k > S_{8,z}$. Hence $p_{z}(k) = \sum_{j=1}^{8} k{8 \choose j} \frac{1}{m^j}e^{\frac{-8k -8 + j + C_{8,z}}{m}}$.\\ Hence $\frac{dp_{z}(k)}{dk} = \left(\sum_{j=1}^{8} {8 \choose j}\frac{1}{m^j}e^{\frac{-8k - 8 + j + C_{8,z}}{m}}\right)\left(1-\frac{8k}{m}\right)$. This is less than $0$ for $k > m/8$. Hence $\forall \ k > m/8$, $p_{z}(k) < p_{z}(m/8)$. Since $p_{z}(m/8)$ is $O(1)$, $p_{z}(k)$ is $O(1)$ for $k > m/8$.\

\textbf{Case 2: } The case where $S_{8,z} > m/8$. First, consider the situation where $S_{1,z} \leq k \leq S_{8,z}$. Let $i$ be such that $S_{i,z} \leq k < S_{i+1,z}$ Now Fact 3 and Remark 5 (see Page 29) imply that $S_{8,z} \leq \frac{6N}{0.06697\sqrt{N}} < 90\sqrt{N}$.
 Hence $k < 90\sqrt{N}$. Therefore, $p_{z}(k) = \sum_{j=1}^{i} {i \choose j} \frac{k}{m^j}e^{\frac{-ik -i + j + C_{i,z}}{m}}$
  $\leq \sum_{j=1}^{i} {i \choose j} \frac{k}{m^j}$
  $\leq \sum_{j = 1}^{i} {i \choose j} \frac{90\sqrt{N}}{m^{j}}$
  $\leq \sum_{j=1}^{i} {i \choose j} \frac{90}{0.06697^j(0.06697\sqrt{2^{30}})^{j-1}}$
  $\leq {8 \choose 1}\frac{90}{0.06697} + \sum_{j=2}^{8} o(1)$, which is $O(1)$.\\ Now when $k > S_{8,z}$,
   $\frac{dp_{z}(k)}{dk} = \left(\sum_{j=1}^{8} {8 \choose j}\frac{1}{m^j}e^{\frac{-8k - 8 + j + C_{8,z}}{m}}\right)(1-\frac{8k}{m})$.
    Since $k > S_{8,z} > m/8$, $\frac{dp_{z}(k)}{dk} < 0 \ \forall \ k > S_{8,z}$. 
   Therefore, $p_{z}(k) < p_{z}(S_{8,z}) \ \forall k \ > S_{8,z}$. Since $p_{z}(S_{8,z})$ is $O(1)$, $p_{z}(k)$ is $O(1) \ \forall \ k > S_{8,z}$. 

 \end{proof}

I will now prove the theorem by approximating the sum of $p_{z}(k)$ over all $k \in \mathbb{N}$ to the integral of $p_{z}(k)$, over each interval $\big[S_{i,z},S_{i+1,z}\big)$. On the interval $[1,S_{1,z})$, $p_{z}(k) = 0$, so $\sum_{k=1}^{S_{i,z}-1} p_{z}(k) = \int_{1}^{S_{1,z}} p_{z}(k)dk$ \textbf{(2)}. Hence for the rest of the proof I will consider $p_{z}(k)$ on intervals $[S_{i,z},S_{i+1,z})$, where $i \geq 1$. Now let $p_{i,z}$ be the function such that $p_{i,z}(k) = k\sum_{j=1}^{i} {i \choose j} \frac{1}{m^j}e^{\frac{-ik -i + j + C_{i,z}}{m}}$ (so $p_{i,z}(k) = p_{z}(k)  \iff k \in [S_{i,z},S_{i+1,z})$). Now by differentiating $p_{i,z}$, we obtain $\frac{dp_{i,z}}{dk} = \left(\sum_{j=1}^i {i \choose j}\frac{1}{m^j}e^{\frac{-ik - i + j + C_{i,z}}{m}}\right)\left(1-\frac{ki}{m}\right)$. Hence $p_{i,z}$ has a single turning point (at $k = \frac{m}{i}$). Since $p_{z}(k) = p_{i,z}(k)$ for some $i$ for every $i \geq 1$, on each interval $[S_{i,z},S_{i+1,z})$, $p_{z}$ is either only decreasing, only increasing, or $\exists$ a $t$ such that $\forall \ k < t$, $p_{z}$ is increasing, and $\forall \ k > t$, $p_{z}$ is decreasing.\\ 
In the case where $p_{z}$ is only decreasing on an interval $[S_{i,z},S_{i+1,z})$, by applying the integral test for convergence, we can conclude that $\int_{S_{i,z}}^{S_{i+1,z}} p_{z}(r)dr < \sum_{k=S_{i,z}}^{S_{i+1,z}-1} p_{z}(k) < \int_{S_{i,z}}^{S_{i+1,z}} p_{z}(r)dr \ + \ \ p_{z}(S_{i,z})$ \textbf{(3)}, while if $p_{z}$ is only increasing on $[S_{i,z},S_{i+1,z})$, $\int_{S_{i,z}}^{S_{i+1,z}} p_{z}(r)dr  -  \ p_{z}(S_{i,z}) < \sum_{k = S_{i,z}}^{S_{i+1,z}-1} p_{z}(k) < \int_{k=S_{i,z}}^{S_{i+1,z}} p_{z}(r)dr$ \textbf{(4)}.\\

Now consider $p_{z}$ on intervals of the form $[S_{i,z},S_{i+1,z})$, where there exists a turning point $t$ such that $\forall \ k < t$, $p_{z}$ is increasing, while $\forall \ k > t$, $p_{z}$ is decreasing (see figure 3.3(a)). Let $j$ be defined such that $j = \max\{m \in \mathbb{N} | n \leq t-1 \}$. Then by applying the integral for convergence test to the intervals $[S_{i,z},j+1]$, and $[j+2,S_{i+1,z}]$, one can see that $\sum_{k = S_{i,z}}^{j+1} p_{z}(k) > \int_{k = S_{i,z}}^{j+1} p_{z}(r)dr$ and $\sum_{k = j+2}^{S_{i+1,z}-1} p_{z}(k) > \int_{j+2}^{S_{i+1,z}} p_{z}(r)dr$. Therefore, $\sum_{S_{i,z}}^{S_{i+1,z}-1} p_{z}(k) + \int_{j+1}^{j+2} p_{z}(r)dr > \int_{S_{i,z}}^{S_{i+1,z}} p_{z}(r)dr$. Now $\int_{j+1}^{j+2} p_{z}(r)dr =$ Ave$\{p_{z}(r) | j+1 \leq r \leq j+2\}$, which is $O(1)$ since all $p_{z}(k)$ are $O(1)$. Hence $\int_{S_{i,z}}^{S_{i+1,z}} p_{z}(r)dr - O(1) < \sum_{S_{i,z}}^{S_{i+1,z}-1} p_{z}(k)$ \textbf{(5)}.\\ Now the integral for convergence test, applied again to the intervals $[S_{i,z},j+1]$, and $[j+2,S_{i+1,z}]$ implies that $\sum_{k = S_{i,z}}^{j} p_{z}(k) < \int_{S_{i,z}}^{j+1} p_{z}(r)dr$, and $\sum_{k = j+3}^{S_{i+1,z}-1} p_{z}(k) < \int_{j+2}^{S_{i+1,z}} p_{z}(r)dr$. Also, it is clear that $p_{z}(j+1) + p_{z}(j+2) < p_{z}(j+1) + p_{z}(j+2) + \int_{j+1}^{j+2} p_{z}(r)dr$. Hence, $\sum_{k = S_{i,z}}^{j} p_{z}(k) + p_{z}(j+1) + p_{z}(j+2) + \sum_{k = j+3}^{S_{i+1,z}-1} p_{z}(k) < \int_{S_{i,z}}^{j+1} p_{z}(r)dr + \int_{j+1}^{j+2} p_{z}(r)dr + \int_{j+2}^{S_{i+1,z}} p_{z}(r)dr$. Therefore, since $p_{z}(j+1)$ and $p_{z}(j+2)$ are $O(1)$, $\sum_{k = S_{i,z}}^{S_{i+1,z}-1} p(k) < \int_{S_{i,z}}^{S_{i+1,z}} p_{z}(r)dr + O(1)$ \textbf{(6)}.

\includegraphics[height = 4.5cm, width = 14cm]{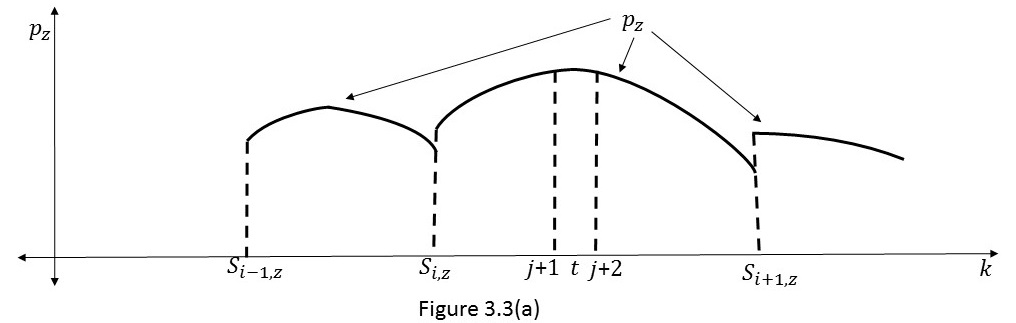}\\

From \textbf{(2)},\textbf{(3)},\textbf{(4)}, \textbf{(5)}, and \textbf{(6)}, we can see that for all intervals $[S_{i,z},S_{i+1,z})$, $$\int_{S_{i,z}}^{S_{i+1,z}} p_{z}(r)dr - O(1) < \sum_{S_{i,z}}^{S_{i+1,z}-1} p_{z}(k) < \int_{S_{i,z}}^{S_{i+1,z}} p_{z}(r)dr + O(1)$$. Therefore, we obtain the required result of $$\int_{1}^{\infty} p_{z}(k)dk - O(1) < \sum_{k=1}^{\infty} p_{z}(k) = \mathrm{E}(Z(z)) <\int_{1}^{\infty} p_{z}(k) dk + O(1)$$ \end{proof}

Hence for large interval sizes $N$, $\int_{1}^{\infty} p_{z}(k) dk$ approximates the expected number of steps until the first collision very well. The following lemma will make the computing of $\int_{1}^{\infty} p_{z}(k) dk$ far more achievable.
\begin{lemma3322} $\forall i,j$ where $1 \leq i \leq 8$, and $j \geq 2$, $\int_{S_{i,z}}^{S_{i+1,z}} k {i \choose j} \frac{1}{m^j} e^{\frac{-ik + C_{i,z} - i + j}{m}}$ is $O(1)$. \end{lemma3322}
\begin{proof}[Proof] The integral of $k {i \choose j} \frac{1}{m^j} e^{\frac{-ik + C_{i,z} - i + j}{m}}$ with respect to $k$ is $$\frac{{i \choose j} (ik + m) e^{\frac{-ik + C_{i,z} - i + j}{m}}}{i^{2}m^{j-1}}$$ so $\int_{S_{i,z}}^{S_{i+1,z}} k {i \choose j} \frac{1}{m^j} e^{\frac{-ik + C_{i,z} - i + j}{m}}dk$ is $${i \choose j}\frac{\left(e^{\frac{-iS_{i,z} + C_{i,z} - i + j}{m}}(iS_{i,z} + m) - e^{\frac{-iS_{i+1,z} + C_{i,z} - i + j}{m}}(iS_{i+1,z} + m)\right)}{m^{j-1}i^{2}} \hspace{4em} \textbf{(7)}$$
Now in the proof of Lemma 3.3.2, I showed that $e^{\frac{-iS_{i,z} + C_{i,z} - i + j}{m}} \leq 1$. Hence \textbf{(7)} is less than $\frac{{i \choose j}(S_{i,z} + m)}{i^{2}m^{j-1}}$.  In the proof of the same lemma , I also showed that $S_{i,z} < 90\sqrt{N}$, and $c_{m} \leq 0.06697\sqrt{N}$. Hence $\frac{{i \choose j}(S_{i,z} + m)}{i^{2}m^{j-1}} < \frac{(90+0.06697)\sqrt{N}}{(0.06697\sqrt{N})^{j-1}}$, which is $O(1)$ for $j \geq 2$. \end{proof}

Therefore, from the formulas at the top of page 33, we see that $$\int_{1}^{\infty} p_{z}(k)dk = \sum_{i=1}^{8} \int_{S_{i,z}}^{S_{i+1,z}} \frac{ik}{m}e^{\frac{-ik + C_{i,z} - i + 1}{m}} + O(1)  \hspace*{12em}\textbf{(8)}.$$ Now each $\int_{S_{i,z}}^{S_{i+1,z}} \frac{ik}{m}e^{\frac{-ik + C_{i,z} - i + 1}{m}}$ is $$e^{\frac{-iS_{i,z} + C_{i,z} - i + 1}{m}}(\frac{m}{i} + S_{i,z}) - e^{\frac{-iS_{i+1,z} + C_{i,z} - i + 1}{m}}(\frac{m}{i} + S_{i+1,z})$$ Substituting this into \textbf{(8)}, we obtain  $$\int_{1}^{\infty} p_{z}(k)dk = \sum_{i=1}^{8} \left(e^{\frac{-iS_{i,z} + C_{i,z} - i + 1}{m}}(\frac{m}{i} + S_{i,z}) - e^{\frac{-iS_{i+1,z} + C_{i,z} - i + 1}{m}}(\frac{m}{i} + S_{i+1,z})\right)   + O(1)$$ 
Now since $C_{i,z} = \sum_{j=1}^{i} S_{j,z}$, $C_{i,z}$ can be expressed as $c_{i,z}\sqrt{N}$, for some $c_{i,z}$ independent of $N$. Hence $$\sum_{i=1}^{8} \left(e^{\frac{-iS_{i,z} + C_{i,z} - i + 1}{m}}(\frac{m}{i} + S_{i,z}) - e^{\frac{-iS_{i+1,z} + C_{i,z} - i + 1}{m}}(\frac{m}{i} + S_{i+1,z})\right)$$ $$= \sum_{i=1}^{8} \left(e^{\frac{(-is_{i,z} + c_{i,z})\sqrt{N} - i + 1}{c_{m}\sqrt{N}}}(\frac{c_{m}}{i} + s_{i,z})\sqrt{N} - e^{\frac{(-is_{i+1,z} + c_{i,z})\sqrt{N} - i + 1}{c_{m}\sqrt{N}}}(\frac{c_{m}}{i} + S_{i+1,z})\right)\sqrt{N}$$ $$= \sum_{i=1}^{8} e^{\frac{-i+1}{m}} \left(e^{\frac{(-is_{i,z} + c_{i,z})}{c_{m}}}(\frac{c_{m}}{i} + s_{i,z}) - e^{\frac{-is_{i+1,z} + c_{i,z}}{c_{m}}}(\frac{c_{m}}{i} + s_{i+1,z})\right)\sqrt{N} \hspace*{4em}\textbf{(9)}$$ Now for the typical interval sizes over which one uses kangaroo methods to solve the IDLP ($N > 2^{30}$), $e^{\frac{i+1}{m}}$ is extremely close to 1. Hence we can safely ignore the $e^{\frac{-i+1}{m}}$ term in \textbf{(9)}. I will state how large the approximation error due to ignoring this $e^{\frac{-i+1}{m}}$ term is when I present my final kangaroo algorithm in section 3.4. Therefore, if we define $c_{z}$ to be $\sum_{i=1}^{8} \left(e^{\frac{(-is_{i,z} + c_{i,z})}{c_{m}}}(\frac{c_{m}}{i} + s_{i,z}) - e^{\frac{-is_{i+1,z} + c_{i,z}}{c_{m}}}(\frac{c_{m}}{i} + s_{i+1,z})\right)$, we have $\int_{1}^{\infty} p_{z}(k)dk - O(1) = c_{z}\sqrt{N}$. Using the result from Theorem 3.3.1.1, we can conclude that \begin{equation} 
\mathrm{E}(Z(z)) = c_{z}\sqrt{N} \pm O(1) \tag{\textbf{10}} 
\end{equation} Hence the expected number of steps until the first collision over all instances  of the IDLP (i.e. over all $z$ with $0 \leq z < N$) is \begin{equation}
\mathrm{E}(Z) = c\sqrt{N} \pm O(1) \tag{\textbf{11}} 
\end{equation} where $c =$ Ave${\left\{ c_{z} | 0 \leq z < 1\right\}}$ $= \int_{0}^{1} c_{z}dz$. Now since at each step, each of the 5 kangaroos make one jump, the expected number of group operations until the first collision occurs across all $z$ is \begin{equation} 5\left(\int_{0}^{1} c_{z}dz\right) \pm O(1) \tag{\textbf{12}}
\end{equation}

\subsubsection{Number of group operations required in Stage 3}
From the analysis provided in [3], if one sets the probability that a group element is distinguished to be $c\log(N)\sqrt{N},$ for some constant $c>0$, the expected number of group operations required in stage 3 is $\sqrt{N}/c\log(N) = o(1)\sqrt{N}$.

\subsubsection{Number of Group operations required in Stage 4}
The kangaroos used in a 2 \textsc{tame}, 2 \textsc{wild1}, and 1 \textsc{wild2} kangaroo method are of the same type as those used in the three kangaroo method. I explained how one may find $z$ from a collision between any of these types of kangaroos when I described the three kangaroo method. In any case, at most 3 addition or subtraction operations modulo $\left|g\right|$, while in the case where a \textsc{Wild1} and a \textsc{Wild2} kangaroo collides, we are required to find $2^{-1} \pmod{\left|g\right|}$. Hence the number of group operations required in this stage is $O(1)$.

Now by summing the number of group operations required in stages 1,2,3 and 4, we obtain the required result that the expected number of group operations required to solve the IDLP by any 2 \textsc{tame}, 2 \textsc{wild1}, and 1 \textsc{wild2} kangaroo method, is approximately $5\left(\left(\int_{0}^1 c_{z} dz\right) + o(1)\right)\sqrt{N} + O(\log(N))$.
\end{proof}

\subsection{Finding a good assignment of starting positions, and average step size}

In this subsection, I will use the formula presented in Theorem 3.3.1 to find the best choice of starting positions and average step size that I could possibly find.\\ The process I will use to do this will be to first state how the formula of Theorem 3.3.1 can be used  compute the running time of an algorithm with particular starting positions and average step size (see Algorithm 1), and then to iterate through various possible starting positions and average step sizes, to find which one minimises the running time.\\ For this purpose, I will define $a$,$b$,$c$,$t_{1}$, and $t_{2}$ to be universal constants independent of the interval size $N$, such that on an interval of size $N$, the 2 \textsc{Wild1} kangaroos start their walks at the positions $aN + zN$ and $cN + zN$, the \textsc{Wild2} kangaroo starts his walk at $bN - zN$, and the 2 \textsc{Tame} kangaroos start their walks at $t_{1}$ and $t_{2}$.  

%\subsubsection*{Computing the running time for an algorithm with particular starting positions and average step size} 
One can see from the formula of Theorem 3.3.1, that finding the running time of any 5 kangaroo method requires finding $\int_{0}^{1} c_{z}dz$. Now calculating $c_{z}$ at any $z$ requires finding $s_{i,z} \ \forall \ i$ with $1 \leq i \leq 8$. The most rigorous approach to finding the average of $c_{z}$ would be to find a formula for each $s_{i,z}$ across all $z$, and to plug this into the formula for $c_{z}$, and then integrate this across all $z$. Finding a direct formula for each $s_{i,z}$ that holds for all possible choices for the starting positions and the average step size proved to be too difficult. I therefore proposed the following simulation based approach for computing $\int_{0}^{1} c_{z}dz$. This is how the \textsc{computeaverage$c_{z}$} function of Algorithm 1 computes $\int_{0}^{1} c_{z}dz$.\\ At a particular $z$, $s_{i,z}$ can be computed by finding the starting positions of all kangaroos on an interval of size $N = 1$ (this means computing $aN + zN = a + z$, $b-z$, $c+z$, $t_{1}$ and $t_{2}$), and then finding the distances between all useful pairs when the kangaroos start at these positions. By then ranking these distances in the manner done in line 5 of Algorithm 1, one can find $D_{i,z}$ for each $i$ between 1 and 8, for an interval of size $N=1$. One can then find $s_{i,z}$ using $D_{i,z}/m = d_{i,z}N/c_{m}\sqrt{N} = s_{i,z}\sqrt{N} = s_{i,z}$. Following this, we can compute all $c_{i,z}$, using $c_{i,z} = \sum_{j=1}^{i} s_{j,z}$. Hence we have all the information we need to compute ${c_{z}}$. By finding the average of $c_{z}$ for a large number of evenly spaced $z$ in $[0,1)$ (so for instance, computing $c_{z}$ for all $z$ in $\left\{z|z = 10^{-p}k,k \in \mathbb{N},0 \leq z < 1\right\}$, with $p \geq 3$), we can approximate $\int_{0}^{1} c_{z}dz$. The pseudocode for the \textsc{matlab\textsuperscript{\textregistered}} function to compute $c_{z}$ is shown in Algorithm 1.

\begin{algorithm}
	\caption{Function for finding the 		average of $c_{z}$}
	\begin{algorithmic}[1]
		\Function{\textsc{computeaverage}$c_{z}$}{$a$,$b$,$c$,$t_{1}$,$t_{2}$,$c_{m}$,$p$}
      		\State $z \longleftarrow 0$
      		\State sumof$c_{z}$ $\longleftarrow 0$
      	
      		\While{$z < 1$}
        		\State distancesarray $\leftarrow$ sort(\par
        		\hskip\algorithmicindent $\big\{\left|(a+z)-t_{1}\right|,\left|(b-z)-t_{1}\right|,\left|(c+z)-t_{1}\right|,\left|(a+z)-t_{2}\right|$\par
        		\hskip\algorithmicindent $\left|(b-z)-t_{2}\right|,\left|(c+z)-t_{2}\right|,\left|(a+z)-(b-z)\right|,\left|(c+z)-(b-z)\right|\big\}$)
        		\For{$i \gets 1 \textrm{ to } 8$}
        		\State $s_{i,z} \leftarrow$ distancesarray$[i]$/$c_{m}$
        		\State $c_{i,z} \leftarrow \sum_{j=1}^{i} s_{j,z}$
        		\EndFor
        		\State $c_{z} \leftarrow \sum_{i=1}^{8} \left(e^{\frac{(-is_{i,z} + c_{i,z})}{c_{m}}}(\frac{c_{m}}{i} + s_{i,z}) - e^{\frac{-is_{i+1,z} + c_{i,z}}{c_{m}}}(\frac{c_{m}}{i} + s_{i+1,z})\right)$
        		\State sumof$c_{z} \leftarrow$ sumof$c_{z} + c_{z}$
        		\State $z \leftarrow z + 10^{-p}$
      		\EndWhile
      		\State \Return{sumof$c_{z}/(10^{p})$}
    	\EndFunction
	\end{algorithmic}
\end{algorithm}

Therefore, if we let $O_{pt}$ denote the output of this function on some specified combination of starting positions and average step size (i.e. values of $a$,$b$,$c$,$t_{1}$,$t_{2}$ and $c_{m}$), then from \textbf{(11)}, the expected number of steps until the first collision occurs for this combination is $O_{pt}\sqrt{N}\pm O(1)$ \textbf{(13)}. Also, from the formula of theorem 3.3.1, the expected number of group operations required to solve the IDLP for our specified combination of starting positions and average step size is approximately $\left(5O_{pt} + o(1)\right)\sqrt{N} \pm O(1)$ \textbf{(14)}. 
I will state how large the error of this approximation is when I present my five kangaroo algorithm in section 3.4.

Therefore, if we let $a_{opt}$,$b_{opt}$,$c_{opt}$,$t_{1_{opt}}$,$t_{2_{opt}}$, and $c_{m_{opt}}$ respectively denote the best values for $a$,$b$,$c$,$t_{1}$,$t_{2}$ and $c_{m}$, then $a_{opt}$,$b_{opt}$,$c_{opt}$,$t_{1_{opt}}$,$t_{2_{opt}}$, and $c_{m_{opt}}$ are the values of $a$,$b$,$c$,$t_{1}$,$t_{2}$ and $c_{m}$ for which the output of the \textsc{computeaverage$c_{z}$} function is smallest.

This fact gives rise the to the following algorithm for finding good values for $a$,$b$,$c$,$t_{1}$,$t_{2}$, and $c_{m}$. The pseudocode for this algorithm is shown in Algorithm 2.
The idea of the algorithm was to start with a range of values for which the optimal values of $a$,$b$,$c$,$t_{1}$,$t_{2}$, and $c_{m}$ lay in. These ranges would be encapsulated in the variables $a_{min}$,$a_{max}$,$b_{min}$,$b_{max}$,$c_{min}$,$c_{max}$,$t_{1_{min}}$,$t_{1_{max}}$,$t_{2_{min}}$ ,$t_{2_{max}}$,$c_{m_{min}}$ and $c_{m_{max}}$, so we would have $a_{min} \leq a_{opt} \leq a_{max}$, $b_{min} \leq b_{opt} \leq b_{max}$, $c_{min} \leq c_{opt} \leq c_{max}$, $t_{1_{min}} \leq t_{1_{opt}} \leq t_{1_{max}}$, $t_{2_{min}} \leq t_{2_{opt}} \leq t_{2_{max}}$. I was unable to prove a range of values for which  $a_{opt}$,$b_{opt}$,$c_{opt}$, $t_{1_{opt}}$,$t_{2_{opt}}$, and $c_{m_{opt}}$ were guaranteed to lie in, but in \textbf{(A)},\textbf{(B)},\textbf{(C)},\textbf{(D)} and \textbf{(E)} (which can be found below), I state and justify some ranges for which $a_{opt}$,$b_{opt}$,$c_{opt}$, $t_{1_{opt}}$,$t_{2_{opt}}$, and $c_{m_{opt}}$ are likely to lie in. Using the \textsc{scanregion} function, I would then find good values for $b$,$c$,$t_{1}$,$t_{2}$ and $c_{m}$ (by \textbf{(A)}, $a$ can be fixed at 0) by computing average$c_{z}$ for evenly spaced (separated by the amount defined by the variable 'gap' in Algorithm 2) values of $b$,$c$,$t_{1}$,$t_{2}$, and $c_{m}$, between the ranges defined by the variables $b_{min}$,$b_{max}$,$c_{min}$,$c_{max}$,$t_{1_{min}}$,$t_{1_{max}}$,$t_{2_{min}}$,$t_{2_{max}}$,$c_{m_{min}}$ and $c_{m_{max}}$ (see lines 3-10). The variables Best$b$,Best$c$, Best$t_{1}$, Best$t_{2}$, and Best$c_{m}$ would represent the values of $b$,$c$,$t_{1}$,$t_{2}$ and $c_{m}$ for which average$c_{z}$ was smallest, across all combinations for which average$c_{z}$ was computed for (this is carried out in lines 9-18).\\ We could then find better values for $b$,$c$,$t_{1}$,$t_{2}$ and $c_{m}$ than Best$b$, Best$c$, Best$t_{1}$, Best$t_{2}$, and Best$c_{m}$, by running the \textsc{scanregion} function on values of $b$,$c$,$t_{1}$,$t_{2}$ and $c_{m}$ in a smaller region centred around Best$b$, Best$c$, Best$t_{1}$, Best$t_{2}$, and Best$c_{m}$ (see lines 25-28), that are separated by a smaller gap (see line 29). By repeating this process multiple times (see lines 24-31), we could keep finding better and better values for $b$,$c$,$t_{1}$,$t_{2}$ and $c_{m}$. Eventually however, the interval for which the \textsc{scanregion} function was called on would shrink to be zero in size. At this point, the the improvements in the best values for $b$,$c$,$t_{1}$,$t_{2}$ and $c_{m}$ between iterations would become negligible. Line 24 determines when this occurs.\\ The algorithm then computes Average$c_{z}$ to a higher degree of accuracy for the optimal values of $b$,$c$,$t_{1}$,$t_{2}$ and $c_{m}$ (see line 32). It does this by setting the variable $p$ in the \textsc{Computeaverage$c_{z}$} function to be 6 ($p$ was set to 3 the main loop (see line 10), due to reasons relating to the practicality of the running time of Algorithm 2).

\begin{algorithm}
	\caption{}
	\begin{algorithmic}[1] 
		\Function{\textsc{scanregion}}{$b_{min}$,$b_{max}$,$c_{min}$,$c_{max}$,$t_{1_{min}}$,$t_{1_{max}}$,$t_{2_{min}}$,$t_{2_{max}}$,$c_{m_{min}}$,$c_{m_{max}}$,gap}
      		
      		\State minaverage$c_{z}$ $\leftarrow \infty$ 
      		\State $b_{values} = \left\{b_{min} + k \times gap | k \in \mathbb{N}, b_{min} \leq b_{min} + k\times gap \leq b_{max}\right\}$
      		\State $c_{values} = \left\{c_{min} + k\times gap | k \in \mathbb{N}, c_{min} \leq c_{min} + k\times gap \leq c_{max}\right\}$
      		\State $t_{1_{values}} = \left\{t_{1_{min}} + k\times gap | k \in \mathbb{N}, t_{1_{min}} \leq t_{1_{min}} + k\times gap \leq t_{1_{max}}\right\}$
      		\State $t_{2_{values}} = \left\{t_{2_{min}} + k\times gap | k \in \mathbb{N}, t_{2_{min}} \leq t_{2_{min}} + k\times gap \leq t_{2_{max}}\right\}$
      		\State $c_{m_{values}} = \left\{c_{m_{min}} + k\times gap | k \in \mathbb{N}, c_{m_{min}} \leq c_{m_{min}} + k\times gap \leq c_{m_{max}}\right\}$
      		\State Combinations = $b_{values} \times c_{values} \times t_{1_{values}} \times t_{2_{values}} \times c_{m_{values}}$
      		\For{each $\left\{b,c,t_{1},t_{2},c_{m}\right\} \in$ Combinations} 
      			\State average$c_{z}$ = \textsc{computeaverage$c_{z}$}(0,$b$,$c$,$t_{1}$,$t_{2}$,$c_{m}$,3)
      			\If{\textsc{average$c_{z}$} $<$ minaverage$c_{z}}$
      				\State minaverage$c_{z} \leftarrow$ average$c_{z}$
      				\State Best$b \leftarrow b$
      				\State Best$c \leftarrow c$
      				\State Best$t_{1} \leftarrow t_{1}$
      				\State Best$t_{2} \leftarrow t_{2}$
      				\State Best$c_{m} \leftarrow c_{m}$
      			\EndIf
      		\EndFor
      		\State \Return minaverage$c_{z}$,Best$b$,Best$c$,Best$t_{1}$,Best$t_{2}$,Best$c_{m}$
      		\EndFunction
      		\State
      		\State $b_{min} \leftarrow -1/2$,$b_{max} \leftarrow 5/2$,$c_{min} \leftarrow 0$,$c_{max} \leftarrow 3$,$t_{1_{min}} \leftarrow t_{2_{min}} \leftarrow -2$,$t_{1_{max}} \leftarrow t_{2_{max}} \leftarrow 9/2$,$c_{m_{min}} \leftarrow 0.06697$,$c_{m_{max}} \leftarrow 0.5330$ 
      		\State gap = $2\times10^{-2}$
      		\State PreviousBestAverage$c_{z} \leftarrow \infty$
      		\State $\big[$CurrentBestAverage$c_{z}$,Best$b$,Best$c$,Best$t_{1}$,Best$t_{2}$,Best$c_{m}\big]$ $\leftarrow$ \textsc{scanregion}($b_{min}$,$b_{max}$,$c_{min}$,$c_{max}$,$t_{1_{min}}$,$t_{1_{max}}$,$t_{2_{min}}$,$t_{2_{max}}$,$c_{m_{min}}$,$c_{m_{max}}$,gap) %\Comment
      		 %SCANREGION($-1/2,5/2,0,3,-2,9/2,-2,9/2,0.06697,0.5330$,gap) \Comment      		
      		\While{PreviousBestAverage$c_{z}$ - CurrentBestAverage$c_{z} \leq 10^{-4}$} 
      		\State $b_{min} \leftarrow$ Best$b$ - gap,	$b_{max} \leftarrow$ Best$b$ + gap, $c_{min} \leftarrow$ Best$c$ - gap,  \State $c_{max} \leftarrow$ Best$c$ + gap, $t_{1_{min}} \leftarrow$ Best$t_{1}$ - gap,  \State $t_{1_{max}} \leftarrow$ Best$t_{1}$ + gap,$t_{2_{min}} \leftarrow$ Best$t_{2}$ - gap,$t_{2_{max}} \leftarrow$ Best$t_{2}$ + gap, \State $c_{m_{min}} \leftarrow$ Best$c_{m}$ - gap,	$c_{m_{max}} \leftarrow$ Best$c_{m}$ + gap 
      		
      		\State gap $\leftarrow$ gap/10 
				\State PreviousBestAverage$c_{z} \leftarrow$ CurrentBestAverage$c_{z}$      			
      			\State $\big[$CurrentBestAverage$c_{z}$,Best$b$,Best$c$,Best$t_{1}$,Best$t_{2}$,Best$c_{m}\big]$ $\longleftarrow$ $\hspace{3em}$\textsc{$\hspace{2em}$scanregion}($b_{min}$,$b_{max}$,$c_{min}$,$c_{max}$,$t_{1_{min}}$,$t_{1_{max}}$,$t_{2_{min}}$,$t_{2_{max}}$,$c_{m_{min}}$,$c_{m_{max}}$,gap)
      		\EndWhile
      		\State AccurateAverage$c_{z}$ = \textsc{computeaverage$c_{z}$}(0,Best$b$,Best$c$,Best$t_{1}$,Best$t_{2}$,Best$c_{m}$,6)
      		\State \Return AccurateAverage$c_{z}$,Best$b$,Best$c$,Best$t_{1}$,Best$t_{2}$,Best$c_{m}$

	\end{algorithmic}
\end{algorithm}

\subsubsection*{Estimates of initial bounds for the optimal values of $a$,$b$,$c$,$t_{1}$,$t_{2}$ and $c_{m}$}

\begin{propA}
$a_{opt} = 0$
\end{propA}
 
\begin{propB}
$\frac{-1}{2} = b_{min} \leq b_{opt} \leq b_{max} = \frac{5}{2}$
\end{propB}

\begin{propC}
$0 = c_{min} \leq c_{opt} \leq c_{max} = 3$, and $c \geq a$
\end{propC}

\begin{propD}
$-2 = t_{1_{min}} = t_{2_{min}} \leq t_{1_{opt}} \leq t_{2_{opt}} \leq t_{1_{max}} = t_{2_{max}} = \frac{9}{2}$
\end{propD}

\begin{propE}
$0.06697 \leq c_{m_{opt}} \leq 0.5330$
\end{propE}

\begin{proof}[Justification of (A)]

Suppose we start our kangaroos walks at $N(d+z)$,$N(b+d-z)$,$N(c+d+z)$,$N(t_{1}+d)$, and $N(t_{2}+d)$, where $d \in \mathbb{R}$. Then for all $z$, one can see that the starting distances between all pairs of kangaroos is the same in this case as when the kangaroos start their walks at $N(0+z)$,$N(b-z)$,$N(c+z)$,$t_{1}N$ and $t_{2}N$ (that is, if we subtract $dN$ from all kangaroos starting positions). Hence if we use the same average step size in both cases, the running time in both algorithms will be the same. Hence for any algorithm where $a > 0$, there exists an algorithm where $a=0$ which has the same running time. Hence we only need to check the case where $a = 0$, so we can claim that $a_{opt} = 0$
\end{proof}

\begin{proof}[Justification of (B)]
The bound presented here is very loose. If $b-a > 2.5$, or $b < a-0.5$ then $\forall \ z$, on an interval of size $N$, the initial distance between the kangaroos that start their walks at $a+z$ and $b-z$ is at least $N/2$ (when $z = 0$, $|(a+z)-(-0.5-z)| = 0.5$, and when $z=1$, $|(a+z)-(2.5-z)| = 0.5$). Now in the three kangaroo method, on an interval of size $N$, the furthest the initial distance between the closest useful pair of kangaroos could be across all $z$ was $N/5$ (this occurred when $z$ was $0,\ 2N/5,\ 3N/5$ and $N$). Now in a good five kangaroo method, since there are more kangaroos (than in a three kangaroo method), we can expect to the furthest initial distance between the closest useful pair of kangaroos across all $z$ to be smaller than $N/5$. Hence if a pair of kangaroos in a five kangaroo method always starts their walks at least distance $N/2$ apart, such a pair is highly unlikely to ever collide before the closest useful pair collides, at any $z$, and will therefore be extremely unlikely to ever be the pair which collides first. In a good 5 kangaroo algorithm, it would be natural to suppose that every useful pair of kangaroos can be the pair whose collision leads to the solving of the IDLP (i.e. can collide first), at some $z$. Hence it is not desirable for there to be a useful pair which always starts its walk at least distance of $0.5N$ apart on an interval of size $N$.
\end{proof}

\begin{proof}[Justification of (C)]
Firstly let $W_{1,1}$ denote the kangaroo that starts its walk at $N(a+z)$, and $W_{1,2}$ denote the kangaroo that starts at $N(c+z)$. Since these kangaroos are of the same types (they're both \textsc{Wild1} kangaroos), if we fix $c$,$t_{1}$,$t_{2}$ and $c_{m}$, then an algorithm where $W{1,1}$ starts at $a+z$ and $W_{1,2}$ starts at $c+z$, has the same running time as an algorithm where $W_{1,1}$ starts at $c+z$ and $W_{1,2}$ starts at $a+z$, since the 2 \textsc{Wild1} walks that start at the same positions in both cases. Hence $c_{opt} \geq a_{opt} = 0$. 
Now $c_{opt} \leq 3$, since if $c > 3$, then on an interval of size $N$, the distance between the kangaroos that start their walks at $N(c+z)$ and $N(b-z)$ can never be less than $N/2$. I gave evidence that this was not desirable in my justification of \textbf{(B)}. Hence we may suppose that $0 \leq c_{opt} \leq 3$.

\end{proof}

\begin{proof}[Justification of (D)]
Firstly, we only need to check values where $t_{2_{opt}} \geq t_{1_{opt}}$, because one can show in a similar way to the one shown in proving $c_{opt} \geq a_{opt}$, that for every case where $t_{2} < t_{1}$, there exists an algorithm with the same running time where $t_{1} < t_{2}$. Now $-2 \leq t_{1},t_{2} \leq 9/2$, since if $t_{1}$ and $t_{2}$ are outside this range, then the starting distance between any wild kangaroo and any of the tame kangaroos on an interval of size $N$ is at least $N/2$. 

\end{proof}

\begin{proof}[Justification of (E)]
I stated in Section 3.2 that by starting the kangaroos walks at $h$, $g^{0.7124N}h$, $g^{1.3562N}h^{-1}$,$g^{0.9274N}$ and $g^{0.785N}$, the lower bound for the expected number of group operations until the closest useful pair collides ($E_{CP}$) of $1.8898\sqrt{N}$ group operations could be realised. Using Lemma 3.2.1 in the proof of Theorem 3.2, the average distance between the closest useful pair across all $z$ when $E_{CP} = 1.8898\sqrt{N}$ group operations is $0.0357N$. In the proof of the same lemma, I also showed that $E_{CP} = 5($Ave$(d(z))/m + m)$, where Ave$(d(z))$ denotes the average distance between the closest useful pair of kangaroos over all $z$. Hence Ave$(d(z)) \geq 0.0357N$ in any 5 kangaroo algorithm. Therefore, $E_{CP} \geq 5(\sfrac{0.0357N}{c_{m}\sqrt{N}} + c_{m}\sqrt{N})$. One could suppose that there must be some bound $B$, such that if the expected number of group operations for the closest useful pair to collide in some algorithm is larger than $B$, then this algorithm could have no chance of being the optimal 5 kangaroo algorithm.  A very loose bound on $B$ is $3\sqrt{N}$. For $E_{CP} \leq 3\sqrt{N}$, we require $5(0.0357N/c_{m}\sqrt{N} + c_{m}\sqrt{N})  \leq 3\sqrt{N}$. For this to occur, $c_{m}$ must range between $0.06697$, and $0.5330$. Hence we can suppose $0.06697 \leq c_{m_{opt}} \leq 0.5330$.

\end{proof}

\subsection*{Remark regarding Lemma 3.3.2} Lemma 3.3.2 required an upper bound on how far apart a useful pair of kangaroos can be when they start their walks. When $t_{2}$ is $4.5$, $b = -0.5$, and $zN = N-1$, the initial distance between the \textsc{Wild2} kangaroo that starts its walk at $(b-z)N$, and the \textsc{Tame} kangaroo that starts his walk at $t_{2}N$, is $6N -1 < 6N$. In any method where the variables $a$, $b$, $c$, $t_{1}$ and $t_{2}$ are within the bounds presented in (A),(B),(C) and (D), no useful pair of kangaroos can start their walks further apart than this.\\ The same lemma also required a lower bound on the average step size used. (E) shows that one lower bound is $0.06697\sqrt{N}$.

\subsection*{Results} When I ran Algorithm 2 in \textsc{matlab\textsuperscript{\textregistered}}, the values for $a$,$b$,$c$,$t_{1}$,$t_{2}$ and $c_{m}$ which were returned had $a = 0$, $b = 1.3578$, $c = 0.7158$, $t_{1} = 0.7930$, $t_{2} = 0.9220$, and $c_{m} = 0.3118$. Algorithm 2 computed average$c_{z}$ to be $0.3474$ in this case.\\ 
Using the result of \textbf{(14)} (see the top of page 40), we see that in an algorithm where $a$,$b$,$c$,$t_{1}$,$t_{2}$ and $c_{m}$ are defined to be these values requires on average $\big(5\times0.3474 + o(1)\big)\sqrt{N} \pm O(1)$ $=\big(1.737 + o(1)\big)\sqrt{N}\pm O(1)$ group operations to solve the IDLP. \\Formula \textbf{(13)} also implies that the expected number of steps until the first collision will be $0.3474\sqrt{N} \pm O(1)$ steps \textbf{(15)}.
%Using the result of \textbf{(17)} \textbf{ADD IN WHAT PAGE THIS APPEARS ON} on, the expected number of group operations required to solve the IDLP on an interval of size $N$, when we start the walks of 5 kangaroos at the group elements $h$, $g^{1.3578N}h^{-1}$, $g^{0.7158N}h$, $g^{0.7930N}$, and $g^{0.422N}$, and use an average step size of $0.3118\sqrt{N}$, is $\left(1.737 + o(1)\right)\sqrt{N}
% \pm O(1)$.\\ If we ignore the approximation error caused by the assumption made in Assumption 3.3, then from the result of \textbf{(16)} on \textbf{Add in page number}, since Algorithm 2 computes average$\overline{c_{z}}$, rather than average$c_{z}$, the multiplicitive factor of lower bound of the approximation error is bounded by a factor of $e^{\frac{-7}{0.3118\sqrt{N}}}$, while the upper bound of the approximation error is bounded by a factor of 1. 
 
\section{Five Kangaroo Method}
I now present my five kangaroo algorithm.\\ 
If we are solving the IDLP on an interval of size $N$ over a group $G$, and we start the walks of five kangaroos at the group elements $h$, $g^{1.3578N}h^{-1}$, $g^{0.7158N}h$, $g^{0.922N}$, and $g^{0.793N}$, and use an average step size of $0.3118\sqrt{N}$, and let the kangaroos walk around the group $G$ in the manner defined in section 3.1, then we can expect to solve the IDLP in on average $\big(1.737 + o(1)\big)\sqrt{N}\pm O(1)$ group operations.\\
Note that $1.3578N$,$0.7158N$,$0.922N$, and $0.793N$ might not be integers, so in practice, the kangaroos would start their walks at $h$, $g^{b}h^{-1}$, $g^{c}h$, $g^{t_{1}}$, and $g_{t_{2}}$, where $b$,$c$,$t_{1}$, and $t_{2}$ are respectively the closest integers to $1.3578N$,$0.7158N$,$0.922N$, and $0.793N$.\\
%The number of group elements that all kangaroos visit before the IDLP is solved is clearly less than the number of group operations required to solve the IDLP. Hence, the number of distinguished points visited by all kangaroos is less than $\big(\big(1.737 + o(1)\big)\sqrt{N}\pm O(1)\big)\c\log(N)/\sqrt{N}$
  %=\big(1.737c + o(1)\big)\log(N)\pm o(1)$. Hence the storage requirements of this algorithm are $O(\log(N)\big)$.\\
There are two main factors which have not been accounted for in this calculation of the running time. The first is eluded to in the remark, stated just before Lemma 3.3.1. This being, that I calculate the expected running time in the instance where for every $z$, the back kangaroo takes the expected number of steps to catch up to the starting position of the front kangaroo in every useful pair of kangaroos, and make the assumption that this is proportional to the average expected running time across all possible walks that the kangaroos can make. I was unable to find a bound for how much the approximation error would be increased by because of this assumption. However, in [3], Galbraith, Pollard and Ruprai make the same assumption in computing the running time of the three and four kangaroo methods. They were also unable to find a bound for much the approximation error would be increased by because of this assumption. However, when Galbraith, Pollard and Ruprai gathered experimental results to test how well their heuristic estimates worked in practice, they found that their experimental results matched their estimated results well [3]. Hence I can assume the magnitude of the approximation error is not affected too badly as a result of this assumption. \\ The other factor that wasn't accounted for in my calculation of the running time, was how I ignore the $e^{\frac{-i+1}{m}}$ term in \textbf{(9)} (see page 37). Since $m = 0.3118\sqrt{N}$ in my 5 kangaroo algorithm, and in practice, one would typically use kangaroo methods on intervals of size at least $2^{30}$, the multiplicative factor of the approximation error due to this is has a lower bound of $e^{\frac{-i+1}{m}} \geq e^{\frac{-7}{0.3118\sqrt{N}}} \geq e^{\frac{7}{0.3118\sqrt{2^{30}}}} \geq 1 - (7\times10^{-4})$, and an upper bound of 1.\\ This five kangaroo algorithm is therefore a huge improvement on the previously optimal five kangaroo method of Galbraith, Pollard, and Ruprai, which required at least $2\sqrt{N}$ group operations to solve the IDLP on average. This algorithm also beats the running time of the three kangaroo method, so it therefore answers one of the main questions of this dissertation. This being, 'Are five kangaroos worse than three?'.  

\chapter{Seven Kangaroo Method}	

Using the same idea that was applied by Galbraith, Pollard and Ruprai in [3] in extending the three kangaroo method to the four kangaroo method, the five kangaroo method can be extended to give a seven kangaroo method, with slightly improved running time.\\
If we are solving the IDLP on an interval of size $N$, let $A$,$B$, and $C$ respectively be the closest \emph{even} integers to $0$,$1.3578N$, and $0.7158N$, and let $T_{1}$ and $T_{2}$ be the closest integers to $0.422N$, and $0.793N$. Suppose we start the walks of 7 kangaroos at the group elements $h$, $g^{B}h^{-1}$, $g^{C}h$, $g^{T_{1}}$, $g^{T_{1}+1}$, $g^{T_{2}}$, and $g^{T_{2}+1}$. Then we are effectively starting two \textsc{Wild1} kangaroos at positions $z$, and $C + z$, one \textsc{Wild2} kangaroo at the position $B - z$, and four \textsc{Tame} kangaroos at the positions $T_{1}$, $T_{1}+1$, $T_{2}$ and $T_{2}+1$. Now $z$, $B-z$ and $C+z$ are either all odd, or all even. Hence all \textsc{Wild1} and \textsc{Wild2} kangaroos start their walks an even distance apart. Also, exactly one of $T_{1}$ and $T_{1}+1$, and exactly one of $T_{2}$ and $T_{2}+1$ are of the same parity as $z$, $B-z$ and $C+z$. I will let $T_{1_{useful}}$ and $T_{2_{useful}}$ denote the \textsc{Tame} kangaroos whose starting positions are of the same parity as the starting positions of the \textsc{Wild1} and \textsc{Wild2} kangaroos, and $T_{1_{useless}}$ and $T_{2_{useless}}$ denote the two other \textsc{tame} kangaroos.\\
Now suppose we make all step sizes even. Then $T_{1_{useless}}$ and $T_{2_{useless}}$ will be unable to collide with any other kangaroo, excluding themselves. However, $T_{1_{useful}}$,$T_{2_{useful}}$, and all \textsc{Wild1} and \textsc{Wild2} kangaroos are able to collide, and are all starting their walks at at almost the exact same positions as the five kangaroos do in the five kangaroo method of section 3.4. Hence, assuming the algorithm is arranged so that all kangaroos jump one after the other in some specified order, then these 5 kangaroos ($T_{1_{useful}}$,$T_{2_{useful}}$, the \textsc{Wild2} kangaroo, and the two \textsc{Wild1} kangaroos) are effectively performing the five kangaroo method of section 3.4, except over an interval of size $N/2$, since the fact that all step sizes are even means that only every second group element is being considered in this method. Hence from the statement of \textbf{(15)} (see the results section on page 44), the expected number of steps until the first collision occurs is $\big(0.3474\sqrt{N/2} \pm O(1)\big)$ $= 0.2456\sqrt{N} \pm O(1)$ steps. However, since there are now 7 kangaroos jumping at each step, the expected number of group operations until the first collision occurs is $7\times0.2456\sqrt{N} \pm O(1)$ = $1.7195\sqrt{N} \pm O(1)$ group operations. By defining the probability that a point is distinguished in the same way as it was in the five kangaroo method, we can conclude that the seven kangaroo method presented here requires on average an estimated $\big(1.7195 + o(1)\big)\sqrt{N} \pm O(1)$ group operations to solve the IDLP.\\ As I have already stated, currently the fastest kangaroo method is the four kangaroo method. This has an estimated average running time of $\left(1.714+o(1)\right)\sqrt{N}$ group operations. Therefore, the seven kangaroo method presented here is very close to being the optimal kangaroo method.

\chapter{Conclusion}

The main results of this thesis were the following. \begin{itemize}
\item \textbf{Result 1:} The presentation of a five kangaroo algorithm which requires on average $\big(1.737 + o(1)\big)\sqrt{N} \pm O(1)$ group operations to solve the IDLP.
\item \textbf{Result 2:} The presentation of a seven kangaroo method that requires on average $\big(1.7195 + o(1)\big)\sqrt{N} \pm O(1)$ group operations to solve the IDLP.
\end{itemize}

For clarity, I will restate the main questions that this thesis attempted to answer here also.

\begin{itemize}
\item \textbf{Question 1:} Can we improve kangaroo methods by using larger numbers of kangaroos?
\item \textbf{Question 2:} Are 5 kangaroos worse than three?
\end{itemize}

Before this dissertation, the four kangaroo method of Galbraith, Pollard, and Ruprai [3] was by the far the best kangaroo algorithm. This algorithm has an estimated average running time of $\big(1.714 + o(1)\big)\sqrt{N}$ group operations. It was unknown whether we could beat the running time of the four kangaroo method by using more than four kangaroos. The fastest algorithm that used more than four kangaroos was far slower than the four kangaroo method, requiring on average at least $2\sqrt{N}$ group operations to solve the IDLP.\\
The five and seven kangaroo algorithms presented in this report are a significant improvement in kangaroo methods that use more than four kangaroos. Even though the running time of the methods presented in this thesis did not beat the running time of the four kangaroo method, they came very close to doing so. Therefore, even though the main question of this dissertation (Question 1) still remains unanswered, we can now have more confidence that kangaroo methods can be improved by using more than four kangaroos.\\ Question 2 was also answered in this thesis, since the estimated average running time of the five kangaroo method presented in section 3.4 ($\big(1.737 + o(1)\big)\sqrt{N} \pm O(1)$ group operations) beat the estimated average running time of the three kangaroo method of [3] ($\big(1.818 + o(1)\big)\sqrt{N}$ group operations).	
% --------------------------------------------------------------

% --------------------------------------------------------------
\renewcommand{\bibname}{References} % changes the header; default: Bibliography

\end{document}